\documentclass{amsart}
\newtheorem{theorem}{Theorem}[section]
\newtheorem{lemma}[theorem]{Lemma}
\newtheorem{proposition}[theorem]{Proposition}

\theoremstyle{definition}
\newtheorem{definition}[theorem]{Definition}
\newtheorem{example}[theorem]{Example}

\theoremstyle{remark}
\newtheorem{remark}[theorem]{Remark}
\numberwithin{equation}{section}
\newcommand{\defi}{\stackrel{\textup{\tiny def}}{=}}
\newcommand{\oo}{\emptyset}

\DeclareMathOperator{\dist}{dist} \DeclareMathOperator{\trace}{trace}
 \DeclareMathOperator{\spa}{span}
\DeclareMathOperator{\ran}{ran} 
\begin{document}
\title[
de Branges -- Rovnyak spaces: the multiscale case] {Schur multipliers and de Branges -- Rovnyak spaces: the
multiscale case}
\author[D. Alpay]{Daniel Alpay}
\address{Department of Mathematics\\
Ben--Gurion University of the Negev\\
Beer-Sheva 84105\\ Israel} \email{dany@math.bgu.ac.il}
\thanks{D. Alpay thanks the Earl Katz family for endowing the chair which
supports his research and NWO, the Netherlands Organization for Scientific Research (grant B 61-524).}
\author[A. Dijksma]{Aad Dijksma}
\address{Department of Mathematics\\
University of Groningen\\
POB 800, NL 9700AV Groningen\\
The Netherlands} \email{dijksma@math.rug.nl}
\thanks{The research of A. Dijksma was in part supported by the Center for
Advanced Studies in Mathematics (CASM) of the Department of Mathematics, Ben--Gurion University}
\author[D. Volok]{Dan Volok}
\address{Department of Mathematics\\
The Weizmann Institute of Science\\
Rehovot 76100\\Israel} \email{danvolok@hotmail.com}

\subjclass{Primary: 93B28;\ 
Secondary: 05C05.}

\keywords{non-commutative power series, system realization, homogeneous tree.}

\begin{abstract}
We consider bounded linear operators acting on the $\ell_2$ space
indexed by the nodes of a homogeneous tree.
Using the Cuntz relations between the primitive shifts on the tree, we generalize the notion of the single-scale
time-varying  point evaluation and introduce the corresponding reproducing
kernel Hilbert space in which
Cauchy's formula holds. These notions are then used in the study of the
Schur multipliers and of the associated
de Branges -- Rovnyak spaces. As an application we obtain realization of Schur multipliers as transfer operators
of multiscale input-state-output systems.
\end{abstract}
\maketitle
\section{Introduction}\label{s1}
In this paper we consider bounded operators acting on the Hilbert space
\begin{equation}\label{l2}{\ell_2(\mathcal{T})}=\{f:{\mathcal{T}}\longrightarrow {\mathbb{C}};\quad
   \|f\|_{\ell_2}^2\defi\sum_{t\in{\mathcal{T}}}|f(t)|^2<\infty\},\end{equation}
where ${\mathcal T}$ is a homogeneous tree of order \(q\ge 1\), that is, an  acyclic, undirected, connected
graph such that every node belongs to exactly \(q+1\) edges (see \cite{MR57:16426}, \cite{MR54:12990}). Such
operators arise in the theory of multiscale linear systems and multiscale stochastic processes. Here we would
like to mention the works \cite{BNW}, \cite{MR1276528}, \cite{MR95f:93057}, where Basseville, Benveniste,
Nikoukhah and Willsky have developed a theory of {\em stationary} multiscale systems and {\em stationary}
multiscale stochastic processes. Connections of their theory with the classical setting when $q=1$ and the tree
$\mathcal{T}$ is the tree of integers $\mathbb{Z}$ (what we shall call the single-scale setting) were explored
in \cite{a-volok} and \cite{MR2182965}. The special case of {\em isotropic} processes was considered in
\cite{MR2196005}; a different approach to isotropic
processes uses the theory of Gelfand pairs (see \cite{MR83k:60019}, \cite{MR563948}).\\

 In what follows we consider the general multiscale setting, without the assumption of stationarity.
 Some of the results presented here were announced in
\cite{adv-cras}.\\

 In order to explain our approach, let us recall that in the stationary single-scale setting
one considers  a function
 \[s(z)=s_0+zs_1+z^2s_2+\dots,\] analytic
and contractive in the open unit disk $\mathbb D$ -- a Schur function. Then the multiplication by $s(z)$ is a
causal contractive operator acting on  the Hardy space $\mathbf{H}_2$ of the unit disk and the kernel
\begin{equation}\label{00ker}K_s(z,w)\defi\dfrac{ 1 - s(z)s(w)^*}
{ 1 - z\overline{w}}\end{equation} is positive in $\mathbb{D}$. The associated reproducing kernel Hilbert space
$\mathbf{H}(s)$ has the form
$$\mathbf{H}(s) =
\sqrt{\mathcal{B}_s}\,\mathbf{H}_2 ;\quad \|\sqrt{\mathcal{B}_s}f\|_{\mathbf{H}(s)} = \|(I -
\pi)f\|_{\mathbf{H}_2},$$ where  $\mathcal{B}_s={I-M_sM_s^*}$ and $\pi$ is
 the orthogonal projection in $\mathbf{H}_2$ onto
$\ker\mathcal{B}_s.$  The space $\mathbf{H}(s)$ is called the de Branges -- Rovnyak space associated with the
Schur function $s$; see \cite{dbr2}, \cite[Appendix]{dbr1}, \cite{sarason94}. It is invariant under the action
of the backward shift operator $R_0,$ defined by \begin{equation}\label{defr0}(R_0f)(z) \defi \dfrac{f(z) -
f(0)}{z}.\end{equation} Moreover, the formulae
\begin{align} \label{real???1}Af &= R_0f,& Bc &= R_0(s\cdot c),\\
\label{real???2} Cf &= f(0),& Dc &= s(0)\cdot c,
\end{align}
where $f\in\mathbf{H}(s)$, $c\in\mathbb{C},$
 define a coisometry
$$\begin{pmatrix} A& B\\ C& D
\end{pmatrix}:
\begin{pmatrix}\mathbf{H}(s)\\ \mathbb{C}\end{pmatrix} \longrightarrow
 \begin{pmatrix}\mathbf{H}(s)\\ \mathbb{C}\end{pmatrix}.$$
In terms of these operators $A,B,C,D$ the Schur function $s(z)$ admits the representation
\begin{equation}
\label{real???} s(z) = D + zC(I - zA)^{-1}B
\end{equation}
and the reproducing kernel $K_s(z,w)$ can be written as
\begin{equation}\label{0ker}
K_s(z,w)=C(I-zA)^{-1}(I-wA)^{-*}C^*.\end{equation} \ \\

Conversely, if ${\mathbf H}$ is a Hilbert space of functions analytic in the open unit disk such that there
exists a coisometry
$$\begin{pmatrix} A& B\\ C& D
\end{pmatrix}:  \begin{pmatrix}\mathbf {H}\\
\mathbb{C}\end{pmatrix} \longrightarrow
 \begin{pmatrix}\mathbf{H}\\ \mathbb{C}\end{pmatrix},$$
 where $A$ and $C$ are as in \eqref{real???1}, \eqref{real???2} (in particular,
 the space $\mathbf{H}$ is $R_0$-invariant),
then the formula \eqref{real???} defines a Schur function $s,$ for which the  kernel $K_s$ is given by
\eqref{0ker} and the associated de Branges -- Rovnyak space coincides with $\mathbf{H}.$  For this and more
general results in the setting of
Pontryagin spaces see \cite[Theorem 3.12 p. 85]{adrs}.\\

The representation \eqref{real???} implies that if a function
\[u(z)=u_0+zu_1+z^2u_2+\dots\in\mathbf{H}_2\]
is given, then the Taylor coefficients $y_0,y_1,y_2,\dots$ of the function
\[y(z)=s(z)\cdot u(z)=y_0+zy_1+z^2y_2+\dots\in\mathbf{H}_2\]
can be recursively determined  as follows:
\begin{equation}\left\{
\begin{aligned}
x_0&=0,\\
x_{n+1}&=Ax_{n}+Bu_n,\\
y_n&=Cx_{n}+Du_n,
\end{aligned}\right.
\label{Pernety}
\end{equation}
where $x_n\in\mathbf{H}(s).$ This fact has many important numerical applications; in the language of system
theory it means that {\em the representation \eqref{real???} is a coisometric realization of the Schur function
$s(z)$ as the transfer function of the input-state-output system \eqref{Pernety}  with the state space
$\mathbf{H}(S)$.}\\

In the  non-stationary single-scale setting, the Hardy space is replaced by the space of upper-triangular
Hilbert -- Schmidt operators, the Schur functions by upper-triangular contractions, the complex variable by the
bilateral shift $Z$ on $\ell_2({\mathbb Z})$ and the constants by diagonal operators (see e.g.
\cite{MR93b:47027}, \cite{MR99g:93001}). In particular, any upper-triangular bounded operator $S$ can be written
as a power series
\begin{equation}\label{sinser}
S= S_{[0]}+ZS_{[1]}+Z^2S_{[2]}+\dots,\end{equation}
 where   $S_{[j]}$ are bounded diagonal operators. In general,
a diagonal operator $D$ does not commute with the shift $Z$. However, they satisfy the commutation relation
$$ZD=D^{(1)}Z,$$
where $D^{(1)}\defi ZDZ^*$ is also a diagonal operator. This fact can be used to define a point evaluation of an
upper-triangular bounded operator at a diagonal ``constant''. At the same time $S$ may be viewed as the input-output operator of a {\em time-varying} causal linear system. In order to construct a non-stationary analogue of the realization \eqref{real???}, it is necessary to consider square-summable sequences of inputs rather than a single input -- in other words, the operator of multiplication by $S$ acting on the space of upper-triangular Hilbert -- Schmidt operators.\\

The multiscale setting considered in this paper can be viewed as the natural multidimensional generalization of
the single-scale non-stationary case. Here expansions of the form \eqref{sinser} are replaced with
non-commutative powers series in $q$ primitive shifts, which satisfy the Cuntz relations (see \cite{MR57:7189},
\cite{MR2003i:42001}). Just as in the single-scale case, the coefficients of these series  do not commute with
the primitive shifts but satisfy certain commutation relations.  Thus the multiscale setting is different from
such multidimensional settings as  the classical theory of formal non-commutative power series with the
coefficients which commute with the indeterminates (see \cite{MR51:583} and  \cite{MR2019348} for recent
developments), the Arveson space of the ball in $\mathbb{C}^n$ (see \cite{MR80c:47010}) and the quaternionic
Arveson space (see  \cite{adr1}, \cite{MR2124899}). In particular, in the last two cases the de Branges --
Rovnyak space associated to a Schur multiplier is Gleason-invariant rather than backward shift-invariant.\\

The paper is organized as follows. Section \ref{s2} is of a review nature. It presents the ordering of the
homogeneous tree $\mathcal{T}$ as introduced by Basseville, Benveniste, Nikoukhah and Willsky and the canonical
representation of a bounded linear operator  on the Hilbert space $\ell_2(\mathcal{T})$  as developed in
\cite{a-volok}. Section \ref{s3} discusses causal operators and, in particular, the algebra of constants. In
Section \ref{s4} we present the point evaluation of a causal operator at a constant. In Section \ref{s5} we
study the space of causal Hilbert -- Schmidt operators which plays here the role of the Hardy space ${\mathbf
H}_2$ of the unit disk. In particular, we present the analogue of Cauchy's formula; see Theorem \ref{proker}.
Schur multipliers, associated kernels, de Branges -- Rovnyak spaces and input-state-output systems are studied
in Section \ref{s6}. In the last section we present the analogue of a Blaschke factor in the
present setting.

\section{Power series representation of bounded operators on $\ell_2({\mathcal T})$}\label{s2}
We start with the ordering of the homogeneous tree $\mathcal{T}$ of order \(q\geq 1.\) Note that, as follows
from the definition (see Introduction), the tree $\mathcal{T}$ is infinite. For each node $t\in\mathcal{T}$ we
consider infinite paths, which begin at $t$. These are infinite sequences of nodes
$$(t_0=t,t_1,t_2,\dots),$$
where each pair of consecutive nodes $t_k,t_{k+1}$ is connected by an edge and each two consecutive edges are
distinct:
$$t_{k+1}\not=t_k\not=t_{k+2},\quad k=0,1,2,\dots.$$
Two such paths
\begin{equation}\label{pa1}(t_0=t,t_1,t_2,\dots)\text{ and
}(s_0=s,s_1,s_2,\dots),\end{equation} which begin at the nodes $t$ and $s$, respectively, are said to be {\em
equivalent} if they coincide modulo finite number of edges: there exist indices $m,n$ such that
\begin{equation}\label{pa2} t_{m+k}=s_{n+k},\quad k=0,1,2,\dots.\end{equation}
 The equivalence classes of paths with
respect to this relation are called the {\em  boundary
points} of the tree \(\mathcal T\).\\

Let us choose and fix some boundary point of  \(\mathcal T\), which will be denoted by\footnote{In the single
scale case, this is $-\infty$.} $\infty_{\mathcal{T}}$. Since the graph \(\mathcal T\) is connected and does not
contain cycles, for each \(t\in\mathcal T\) there exists a unique representative of the equivalence class
$\infty_{\mathcal{T}}$, which begins at the node
\(t\).\\

For a pair of nodes \(t,s\), let the corresponding representatives of the boundary point $\infty_{\mathcal{T}}$
be given by \eqref{pa1}. They
 coincide modulo finite number of edges, that is, \eqref{pa2} holds for
some $m$ and $n$. Let us choose the minimal $m$ and $n$ for which \eqref{pa2} holds. Then we denote the node
$t_m=s_n$ by $ s\wedge t$ and call the number $m+n$ {\em the distance} \(\dist(s,t)\) between the nodes
\(s\) and \(t\).\\

Using these notations, we define the  partial ordering $\preceq$ and the equivalence relation $\asymp$ as
follows:
\begin{align}\label{defparord}
s&\preceq t, &\text{if}\quad\dist(s,s\wedge t)&\leq \dist(t,s\wedge t).
\\ \label{defeqrel}
s&\asymp t, &\text{if} \quad\dist(s,s\wedge t)&= \dist(t,s\wedge t).
\end{align}
The equivalence classes with respect to the relation $\asymp$ are called {\em horocycles.}\\

Furthermore, we choose and fix $q$ mappings
\[\alpha_1,\ldots,\alpha_q:\mathcal{T}\longrightarrow\mathcal{T}
\]
acting {\em on the right}
\[
\{t\alpha_1,\dots,t\alpha_q\}=\{s\in\mathcal{T}\,:\,t\preceq s,\ \dist(t,s)=1\}.\] The mappings
$\alpha_1,\ldots,\alpha_q$ are called {\em primitive shifts}.\\

   The induced {\em left} action of the primitive shifts \(\alpha_1,\dots,\alpha_q\)
   on the space\footnote{See \eqref{l2}.}
   ${\ell_2(\mathcal{T})}$ is given by:
   \begin{equation}\label{prsh}\big({\alpha_j}f\big)(t)\defi f(t\alpha_j
),\quad f\in\ell_2(\mathcal{T}),\ t\in\mathcal{T},\ j=1,\dots,q.\end{equation} Thus the primitive shifts
\(\alpha_1,\dots,\alpha_q\) can be also viewed as
   bounded linear operators on $\ell_2(\mathcal{T})$.
    The adjoint operators are given
by
\begin{equation}
\label{prshadj}\big(\alpha_j^*f\big)(t)=\left\{\begin{aligned}
f(s),&\quad\text{if }t\text{ is of the form } t=s\alpha_j,\\
0,&\quad\text{otherwise}.\end{aligned} \right.\end{equation} and satisfy the Cuntz relations:
\begin{equation}\label{cuntz}
{\alpha_i}\alpha_j^*=\delta_{i,j}\cdot I,\quad \sum_{j=1}^q\alpha_j^*\alpha_j=I.\end{equation} In other words,
the operator matrix
\[{\alpha}\defi\begin{pmatrix}\alpha_1\\ \alpha_2\\
\vdots \\ \alpha_q
\end{pmatrix}: \ell_2(\mathcal{T})\longrightarrow
\ell_2(\mathcal{T})^q\] is unitary:
\begin{equation}\label{comcuntz}
{\alpha}{\alpha}^*=I_{\ell_2({\mathcal T})^q},\quad {\alpha}^*{\alpha}=I_{\ell_2({\mathcal T})}.
\end{equation}

Since \(\mathcal{T}\) is a tree, the shifts \(\alpha_j\) form a free semigroup, which we  denote by
\(\mathcal{F}_q\). Every  element $w\in\mathcal{F}_q$ acts on the tree $\mathcal{T}$ on the right: $$t\mapsto
tw,$$ and on the space $\ell_2(\mathcal{T})$ on the left: $$f\mapsto wf,\quad (wf)(t)=f(tw).$$ The unit element
of $\mathcal{F}_q$ will be denoted by $\oo.$ For $w\in\mathcal{F}_q$ we  also use the notation:
\begin{equation}
|w|\defi\left\{\begin{aligned} 0,&\quad\text{if }w=\oo,\\
n,&\quad\text{if }w=\alpha_{i_1}\dots\alpha_{i_n}.
\end{aligned}
\right.
\end{equation}

\begin{definition}\label{reduce}
A pair of elements $w_1,w_2\in\mathcal{F}_q$ is said to be {\em reducible} if $w_1$ and $w_2$ can be represented
as
\[w_1=\alpha_iv_1,\quad w_2=\alpha_iv_2\]
for some $v_1,v_2\in\mathcal{F}_q$ and some primitive shift $\alpha_i$.\end{definition}

\begin{remark}\label{irredu}
Note that a pair of elements $w_1,w_2\in\mathcal{F}_q$ is irreducible if and only if there exists
$t\in\mathcal{T}$ such that
\begin{equation}\label{triv56}
(tw_1)\wedge(tw_2)=t.\end{equation} In this case \eqref{triv56} holds for all $t\in\mathcal{T}.$
\end{remark}

 Let $\mathbf{X}(\mathcal{T})$ denote the $C^*$-algebra of
   bounded linear operators on $\ell_2(\mathcal{T})$.
The elements of the semigroup $\mathcal{F}_q$  appear in the non-commutative power series representations of
elements of $\mathbf{X}(\mathcal{T})$. The coefficients of these power series are diagonal operators with
respect to the standard basis $\chi_t$ of
$\ell_2(\mathcal{T}): $\[\chi_t(s)\defi \left\{\begin{aligned} 1,&\quad\text{if }t=s,\\
0,&\quad\text{otherwise}.\end{aligned} \right.\]  The precise result can be formulated as follows:
\begin{theorem}\label{funthrm}
Every  operator \(S\in\mathbf{X}(\mathcal T)\) can be represented in the form
\begin{equation}\label{canrep}
S=\sum_{w_1,w_2\in\mathcal{F}_q }^\prime w_1^*{w_2}S_{{w_1},{w_2}},
\end{equation}
where:
\begin{enumerate}
\item
 \( S_{w_1,w_2}\)  are elements  of $\mathbf{X}(\mathcal{T})$ which are diagonal with respect to the
standard basis $\chi_t$ of $\ell_2(\mathcal{T}).$

\item The notation $$\sum_{w_1,w_2\in\mathcal{F}_q }^\prime$$ means that
the summation is taken over all the irreducible pairs $w_1,w_2\in\mathcal{F}_q.$

\item Convergence is {\em absolute pointwise}:
for every $f\in{\ell_2(\mathcal{T})}$ and $t\in\mathcal{T}$ the series
\[\sum_{{w_1},{w_2}\in\mathcal{F}_q}^\prime
\left({w_1}^*{w_2}S_{{w_1},{w_2}}f\right)(t),\]  is absolutely convergent and its sum is equal to $(Sf)(t)$.
\end{enumerate}

Moreover, the diagonal coefficients \( S_{w_1,w_2}\) of the series \eqref{canrep} are partially determined by
\begin{equation}\label{undet}
(S_{w_1,w_2}\chi_{r})(r)= ( S\chi_{r})(sw_1) \quad\text{if }r\in\mathcal{T}\text{ is of the form } r=sw_2
\end{equation} The rest of the diagonal entries of \( S_{w_1,w_2}\), that is, the values
$(S_{w_1,w_2}\chi_{r})(r)$ for $r\not\in\mathcal{T}w_2$ can be assigned arbitrarily (as long as they are
uniformly bounded).
\end{theorem}
\begin{proof}
Let $S\in\mathbf{X}(\mathcal{T}),$ \(f\in\ell_2(\mathcal{T})\) and \( t\in\mathcal{T}\) be fixed. Then
\[(Sf)(t)=\sum_{r\in\mathcal{T}}(
S\chi_r)(t)\cdot f(r),\] where the series is absolutely convergent because of the Cauchy -- Schwarz inequality.

Observe that, in view of Remark \ref{irredu}, for each $r\in\mathcal{T}$ there exists a unique triple
$s\in\mathcal{T},w_1,w_2\in\mathcal{F}_q$ such that
\[t=sw_1,\quad r=sw_2\]
and the pair $w_1,w_2$ is irreducible. (In particular, $s=t\wedge r.$) Therefore, we have
\[
(Sf)(t)=\sum_{\substack{s\in\mathcal{T},w_1,w_2\in\mathcal{F}_q:\\t=sw_1}}^\prime (S\chi_{sw_2})(sw_1)\cdot
f(sw_2).\] Let now  \( S_{w_1,w_2}\)  be elements  of $\mathbf{X}(\mathcal{T})$ which are diagonal with respect
to the standard basis $\chi_t$ of $\ell_2(\mathcal{T})$ and satisfy \eqref{undet}. Then one can rewrite the last
identity as
\[
(Sf)(t) =\sum_{\substack{s\in\mathcal{T},w_1,w_2\in\mathcal{F}_q:\\t=sw_1}}^\prime
 (S_{w_1,w_2}f)(sw_2)\\=\sum_{\substack{s\in\mathcal{T},w_1,w_2\in\mathcal{F}_q:\\t=sw_1}}^\prime
 (w_1^*w_2S_{w_1,w_2}f)(t).\]
But, in view of \eqref{prshadj}, for every $w\in\mathcal{F}_q$ and $g\in\ell_2(\mathcal{T})$ we have
\begin{equation}\label{adj}
(w^*g)(t)=\left\{\begin{aligned} g(s),&\quad\text{if }t\text{ is of the form } t=sw,\\
0,&\quad\text{otherwise}.\end{aligned} \right.\end{equation}
Hence  the identity \eqref{canrep} holds in the  sense of pointwise absolute  convergence.\\

Furthermore, let \(t,s\in\mathcal T\) be fixed and let \(S\in\mathbf X(\mathcal T)\) admit (in the  sense of the
pointwise absolute  convergence) a representation of the form \eqref{canrep}, where the coefficients
\(S_{w_1,w_2}\in\mathbf X(\mathcal T)\) are diagonal with respect to the standard basis $\chi_t$ of
$\ell_2(\mathcal{T}).$ Then
\[
(S\chi_t)(s)=\sum_{{w_1},{w_2}\in\mathcal{F}_2}^\prime \left({w_1}^*{w_2}S_{{w_1},{w_2}}\chi_t\right)(s),
\]
but, in view of \eqref{adj}, the terms of the series on the right-hand side satisfy the relations
\[\left({w_1}^*{w_2}S_{{w_1},{w_2}}\chi_t\right)(s)=
\left\{\begin{aligned} (S_{{w_1},{w_2}}\chi_t)(t),&\quad\text{if }t=(t\wedge s)w_2\text{ and }
s=(t\wedge s)w_1,\\
0,&\quad\text{otherwise},\end{aligned} \right.\] hence the series contains at most one non-zero term and
\eqref{undet}
 follows.
\end{proof}
\begin{remark}\label{unrem}
In order to avoid ambiguity, we shall usually normalize the diagonal coefficients $S_{w_1,w_2}$ of the
representation \eqref{canrep} for an operator $S\in\mathbf{X}(\mathcal{T})$ as follows:
 \begin{equation}\label{unass}S_{w_1,w_2}\chi_r=0,\quad\text{if }r\not\in\mathcal{T}w_2.
\end{equation}
\end{remark}

Since the diagonal coefficients \(S_{w_1,w_2}\) in the expansion \eqref{canrep} do not commute, in general, with
the shift operators $w_1^*w_2$, it is of interest to study the representations \eqref{canrep} in the special
cases when $S$ is of the form $S=Dw^*$ or $S=Dw,$ where $D\in\mathbf{X}(\mathcal{T})$ is a diagonal operator
with respect to the standard basis $\chi_t$ of $\ell_2(\mathcal{T})$ and $w\in\mathcal{F}_q$.

\begin{lemma}\label{trlem1}
Let $D\in\mathbf{X}(\mathcal{T})$ be a diagonal operator with respect to the standard basis $\chi_t$ of
$\ell_2(\mathcal{T})$ and let $w\in\mathcal{F}_q$. Then
\[Dw^*=w^*D^\prime,\quad Dw=wD^{\prime\prime},\]
where  $D^\prime,D^{\prime\prime}\in\mathbf{X}(\mathcal{T})$ are diagonal operators given by
\begin{align*}
D^\prime\chi_t&=(
D\chi_{tw})(tw)\cdot\chi_t\quad\forall t\in\mathcal{T}\\
D^{\prime\prime}\chi_t&= \left\{\begin{aligned} (D\chi_{s})(s)\cdot\chi_t,&\quad\text{if }t\in\mathcal{T}
\text{ is of the form } t=sw,\\
0,&\quad\text{otherwise}.\end{aligned} \right.\end{align*}
\end{lemma}
\begin{proof}
 In view of \eqref{prsh}, \eqref{prshadj} we have
\begin{equation}\label{fqw}w^*\chi_t=\chi_{tw},
\quad w\chi_t=\left\{\begin{aligned} \chi_s,&\quad\text{if }t\text{ is of the form } t=sw,\\
0,&\quad\text{otherwise}.\end{aligned} \right.\end{equation} Since the operators $D,D^\prime, D^{\prime\prime}$
are diagonal, the rest of the proof is straightforward.
\end{proof}

\section{Causal bounded operators and constants}\label{s3}
\begin{definition}\label{acau} Let $S\in\mathbf{X}(\mathcal{T}).$
$S$ is said to be causal if for every node $s\in\mathcal{T}$ and every element $f\in\ell_2(\mathcal{T})$ such
that $$f(t)=0\quad\forall t\preceq s$$ it holds that $$(Sf)(t)=0\quad\forall t\preceq s.$$
\end{definition}

\begin{example} For every $w\in\mathcal{F}_q$ the adjoint operator
$w^*\in\mathbf{X}(\mathcal{T})$ is causal, as follows from \eqref{fqw}.
\end{example}

\begin{proposition}\label{causal}
Let $S\in\mathbf{X}(\mathcal{T})$ be represented by the pointwise absolutely convergent series
\[
S=\sum_{{w_1},{w_2}\in\mathcal{F}_q}^\prime w_1^*{w_2}S_{{w_1},{w_2}},\] where \(S_{w_1,w_2}\in\mathbf
X(\mathcal T)\) are diagonal operators with respect to the standard basis $\chi_t$ of $\ell_2(\mathcal{T}),$
normalized by \eqref{unass}.

 Then $S$ is causal if and only if
$$S_{w_1,w_2}=0\quad \text{whenever}\quad |w_1|<|w_2|.$$
\end{proposition}
\begin{proof}
First let us assume that $S$ is causal. Let $w_1,w_2\in\mathcal{F}_q$ be an irreducible pair such that
$|w_1|<|w_2|$ and let $s\in\mathcal{T}$. Then, according to \eqref{defparord} and \eqref{defeqrel},
$$sw_1\preceq sw_2\quad\text{and}\quad sw_1\not\asymp sw_2.$$ Therefore, from Definition \ref{acau} and
 the formula \eqref{undet}
of Theorem \ref{funthrm} it follows that
\[S_{w_1,w_2}\chi_{sw_2}=
(S\chi_{sw_2})(sw_1)\cdot\chi_{sw_2}=0.\] In view of \eqref{unass}, we conclude that $S_{w_1,w_2}=0.$

Conversely, assume that $S_{w_1,w_2}=0$ whenever $|w_1|<|w_2|.$ Let $s,t\in\mathcal{T}$ be such that $t\preceq
s$ and $t\not\asymp s$. Then there exists a unique  pair of elements $w_1,w_2\in\mathcal{F}_q$ such that
\[t=(t\wedge s)w_1,\quad s=(t\wedge s)w_2.\] By definition of $t\wedge s$, this pair $w_1,w_2$ is irreducible.
In view of  \eqref{defparord} and \eqref{defeqrel}, $|w_1|<|w_2|$. Hence, according to the formula \eqref{undet}
of Theorem \ref{funthrm},
\[(S\chi_s)(t)=(S_{w_1,w_2}\chi_s)(s)=0.\]

Thus $(S\chi_s)(t)=0$ for every pair of nodes $s,t\in\mathcal{T}$  such that $t\preceq s$ and $t\not\asymp s$.
In view of Definition \ref{acau}, this means that $S$ is causal.
\end{proof}

 We shall denote the algebra of causal operators $S\in\mathbf{X}(\mathcal{T})$ by ${\bf H}({\mathcal T})$. Note that,
 in
view of Definition \ref{acau}, the algebra $\mathbf{H}(\mathcal{T})$ is  closed  in $\mathbf{X}(\mathcal{T})$ in
the pointwise sense: if a sequence $S_1,S_2,\dots$ of elements of $\mathbf{H}(\mathcal{T})$ and an element
$S\in\mathbf{X}(\mathcal{T})$ are such that for every $f\in\ell_2(\mathcal{T})$ and $t\in\mathcal{T}$
$\lim_{n\rightarrow \infty}(S_nf)(t)=(Sf)(t),$ then $S\in\mathbf{H}(\mathcal{T}).$

In order to study the algebra $\mathbf{H}(\mathcal{T})$ further, we consider its subalgebra
\[\mathcal{C}\defi\{S\in\mathbf{X}(\mathcal{T})\, : \,
S,S^*\in\mathbf{H}(\mathcal{T})\}.\] Elements of $\mathcal{C}$ play the role of constants in the present
setting; we note that
\begin{equation}
\label{const}S\in\mathcal{C}\quad\Leftrightarrow\quad S\chi_t\in\overline{\spa}\{\chi_s\,:\,s\asymp t\}\ \forall
t\in\mathcal{T},\end{equation} where $\overline{\spa}$ denotes closed linear span. Thus the subalgebra
$\mathcal{C}$ is  closed  in $\mathbf{X}(\mathcal{T})$ in the pointwise sense.

\begin{remark}\label{charac}
Note that,  according to Proposition \ref{causal} and Theorem \ref{funthrm}, an element $S$ of the algebra
$\mathbf{X}(\mathcal{T})$ belongs to the subalgebra $\mathcal{C}$
 if and only if it is of the form
\[S=\sum_{\substack{{w_1},{w_2}\in\mathcal{F}_q\\|w_1|=|w_2|}}^\prime
w_1^*{w_2}S_{{w_1},{w_2}},\] where $S_{{w_1},{w_2}}$ are diagonal operators with respect to the standard basis
$\chi_t$ of $\ell_2(\mathcal{T}).$
\end{remark}

\begin{theorem}\label{trivprop}
Let \(S\in\mathbf{H}({\mathcal T}).\)  Then \(S\) can be represented as the  pointwise absolutely convergent
series
\begin{equation}
\label{llg} S=\sum_{w\in\mathcal{F}_q}w^*S_{[w]},\end{equation} where \(S_{[w]}\in{\mathcal{C}}\) are uniquely
determined by
\begin{equation}\label{tayco}
(S_{[w]}\chi_t)(s)=\left\{\begin{aligned}(
S\chi_t)(sw),&\quad\text{if}\quad t\asymp s,\\
0,&\quad\text{otherwise},\end{aligned}\right.\end{equation} and satisfy
\begin{equation}\label{rough}
\|S_{[w]}\|\leq\left\|\sum_{\substack{v\in\mathcal{F}_q\\|v|= |w|}}v^*S_{[v]}\right\|\leq\|S\|.
\end{equation}
\end{theorem}
In the proof of Theorem \ref{trivprop} we shall use the following lemma:
\begin{lemma}\label{coco}
Let $T\in\mathbf{X}(\mathcal{T})$  and let
\begin{equation}\label{llg3}
T=\sum_{{w_1},{w_2}\in\mathcal{F}_q}^\prime w_1^*{w_2}T_{{w_1},{w_2}},\quad\text{where}\quad
T_{{w_1},{w_2}}\text{ are diagonal,}\end{equation} be the pointwise absolutely convergent expansion of $T$ as in
Theorem \ref{funthrm}. Then the series
\[T_{[\oo]}\defi\sum_{\substack{{w_1},{w_2}\in\mathcal{F}_q\\|w_1|=|w_2|}}^\prime
w_1^*{w_2}T_{{w_1},{w_2}}\] converges pointwise absolutely  in  $\mathcal{C}$ and
\begin{equation}\label{inecoco}\|T_{[\oo]}\|\leq\|T\|.\end{equation}
\end{lemma}
\begin{proof}
First we observe that, since the series \eqref{llg3} is pointwise absolutely convergent, the series
\[(T_{[\oo]}f)(t)\defi\sum_{\substack{{w_1},{w_2}\in\mathcal{F}_q\\|w_1|=|w_2|}}^\prime
(w_1^*{w_2}T_{{w_1},{w_2}}f)(t)\] converges absolutely for each $f\in\ell_2(\mathcal{T})$ and
$t\in\mathcal{T}$.\\

 Now let
$$f=\sum_{i=1}^n f_i\chi_{t_i},\quad\text{where}\quad f_i\in\mathbb{C},$$
 let $\mathfrak{h}_1,\dots \mathfrak{h}_k$ be horocycles such that
$$\{t_1,\dots,t_n\}\subset \mathfrak{h}_1\cup \dots\cup \mathfrak{h}_k$$ and let $\pi_1,\dots,\pi_k$ denote the
corresponding  orthogonal projections in $\ell_2(\mathcal{T})$:
\[\pi_j=\sum_{t\in
\mathfrak{h}_j}\langle\cdot,\chi_t\rangle_{\ell_2(\mathcal{T})}\,\chi_t.\] For each $j$  the relations
\eqref{fqw} imply
\[(T_{[\oo]}\pi_jf)(t)=\left\{\begin{aligned}
(\pi_j T\pi_jf)(t),&\quad\text{if}\quad t\in\mathfrak{h}_j,\\
0,&\quad\text{otherwise},\end{aligned}\right.
\] hence
\[T_{[\oo]}\pi_jf=\pi_jT\pi_jf\in\ran(\pi_j).\]
Since \[f=\sum_{j=1}^k\pi_jf,\] we observe that \(T_{[\oo]}f\in\ell_2(\mathcal{T})\) and, moreover,
\begin{multline*}
\|T_{[\oo]}f\|_{\ell_2(\mathcal{T})}^2=\left\|\sum_{j=1}^kT_{[\oo]}\pi_jf\right\|_{\ell_2(\mathcal{T})}^2=
\sum_{j=1}^k\|\pi_jT\pi_jf\|_{\ell_2(\mathcal{T})}^2\\\leq\sum_{j=1}^k\|T\|^2\,\|\pi_jf\|_{\ell_2(\mathcal{T})}^2=
\|T\|^2\,\|f\|_{\ell_2(\mathcal{T})}^2.\end{multline*} Since $\spa\{\chi_t\,:\,t\in\mathcal{T}\}$ is dense in
\(\ell_2(\mathcal{T})\), we conclude that $$T_{[\oo]}\in\mathbf{X}(\mathcal{T}),\quad \|T_{[\oo]}\|\leq\|T\|.$$
Finally, in view of Remark \ref{charac}, $T_{[\oo]}\in\mathcal{C}.$
\end{proof}

\begin{proof}[Proof of Theorem \ref{trivprop}]
According to Theorem \ref{funthrm} and Proposition \ref{causal}, $S$ admits in the sense of  pointwise absolute
convergence  a representation of the form
\[S=\sum_{\substack{{w_1},{w_2}\in\mathcal{F}_q\\|w_1|\geq|w_2|}}^\prime w_1^*{w_2}S_{{w_1},{w_2}},\quad\text{where}\quad
S_{{w_1},{w_2}}\text{ are diagonal}.\]
 Denote
$$S_{[w]}\defi (wS)_{[\oo]}, \quad w\in\mathcal{F}_q.$$ Then,
 as follows from \eqref{cuntz}, $S_{[w]}$ admits in the sense of pointwise absolute  convergence the
representation
\[S_{[w]}=\sum_{\substack{{w_1},{w_2}\in\mathcal{F}_q\\|w_1|=|w_2|}}^\prime
w_1^*{w_2}S_{{w_1w},{w_2}}.\] Hence $S$ admits the representation \eqref{llg} and the relations \eqref{undet} in
Theorem \ref{funthrm} imply \eqref{tayco}.\\

 Finally, we note that, as follows from Lemma \ref{trlem1} and
 \eqref{cuntz},
\[S_{[w]}=w\sum_{\substack{v\in\mathcal{F}_q\\|v|= |w|}}v^*S_{[v]}\quad\text{and}\quad(Sw)_{[\oo]}=\sum_{\substack{v\in\mathcal{F}_q\\|v|=
|w|}}v^*S_{[v]}w .\] Using the identity  $ww^*=I$ and the inequality \eqref{inecoco} in Lemma \ref{coco}, we
obtain
\[\|S_{[w]}\|
\leq\left\|\sum_{\substack{v\in\mathcal{F}_q\\|v|= |w|}} v^*S_{[v]}\right\| =\|(Sw)_{[\oo]}w^*\| \leq
\|Sw\|\leq\|S\|.\] Thus the inequality \eqref{rough} holds.
\end{proof}

\begin{remark}\label{linexp1} Let $c\in\mathcal{C}$, let $S,T\in\mathbf{H}(\mathcal{T})$ and let
\[S=\sum_{w\in\mathcal{F}_q}w^*S_{[w]},\quad T=\sum_{w\in\mathcal{F}_q}w^*T_{[w]},\quad\text{where}\quad
S_{[w]},T_{[w]}\in\mathcal{C},\] be the pointwise absolutely convergent expansions of $S,T$ as in Theorem
\ref{trivprop}.  Then
\[Sc+T=\sum_{w\in\mathcal{F}_q}w^*(S_{[w]}c+T_{[w]}),\]
where convergence is again pointwise absolute.
\end{remark}

Note that the  coefficients \(S_{[w]}\in\mathcal{C}\) in the expansion \eqref{llg} do not commute, in general, with
the shift operators $w^*$. The following lemma deals with the special case when
$S$ is of the form $S=Cw^*$, where $C\in\mathcal{C}.$

\begin{lemma}\label{trivlem2}
Let \(w\in\mathcal{F}_q\)  and let \(C\in\mathcal{C}\). Then
\begin{equation}\label{stnot1} Cw^*=\sum_{\substack{v\in\mathcal{F}_q\\|v| =|w|}}{v}^*
C_v,\quad\text{where}\quad C_v= v C w^*\in\mathcal{C}.\end{equation}
\end{lemma}

\begin{proof}
Let  $v\in\mathcal{F}_q$ be such that $|v|=|w|$ and consider the operator $C_v=vCw^*\in\mathbf{X}(\mathcal{T})$.
For each
 $t\in\mathcal{T}$
the relations \eqref{const} and \eqref{fqw} imply that
\begin{multline*}C_v\chi_t=vC\chi_{tw}\in\overline{\spa}\{v\chi_u\,:\,u\asymp tw\}\\
=\overline{\spa}\{\chi_s\,:\,sv\asymp tw\}\subset \overline{\spa}\{\chi_s\,:\,s\asymp t\}\end{multline*} and
hence, according to \eqref{const}, $C_v\in\mathcal{C}.$\\

 Furthermore, we note that the Cuntz relations \eqref{cuntz} imply
\begin{equation}\label{trivid}
\sum_{\substack{v\in\mathcal{F}_q\\|v| =n}}v^*v =I,\quad n=0,1,2,\dots\end{equation} Hence
\[Cw^*=\sum_{\substack{v\in\mathcal{F}_q\\|v|
=|w|}}v^* v Cw^*\] and \eqref{stnot1} follows.
\end{proof}

\section{The point evaluation of causal operators}\label{s4}
In this section we  define a point evaluation for the elements of
$\mathbf{H}(\mathcal{T})$ at the "points" from $\mathcal{C}.$\\

We consider the set of $q$-tuples
\begin{equation}\label{ball}
\mathbb{B}(\mathcal{T})\defi\{c =\begin{pmatrix}{c}_1&\dots&{c}_q\end{pmatrix}\in{\mathcal
C}^{q}\,:\,\lim_{n\rightarrow\infty}\|(c\alpha)^n\|^{\frac{1}{n}}<1\},\end{equation} which plays the role of the
unit disk in the present setting.\\

Let ${c}\in{\mathbb{B}}({\mathcal{T}})$ and let $S\in\mathbf{H}(\mathcal{T})$. Let
\begin{equation}\label{llg2}S=\sum_{w\in\mathcal{F}_q}w^*S_{[w]},
\quad\text{where}\quad S_{[w]}\in\mathcal{C},\end{equation}
 be the pointwise absolutely convergent expansion  of
$S$ as in Theorem \ref{trivprop}.
 Then, in view of the estimate \eqref{rough},
 the series
 \begin{equation}\label{lepoev}
 S^\wedge(c)\defi\sum_{n=0}^\infty(c\alpha)^n\left(\sum_{\substack{w\in\mathcal{F}_q\\|w|=n}}w^*S_{[w]}\right)
 \end{equation}
 converges
absolutely with respect to the operator norm. It follows from Lemma \ref{trivlem2} and the Cuntz relations
\eqref{cuntz} that each term of the series \eqref{lepoev} belongs   to the algebra of constants $\mathcal{C},$
which is
closed in $\mathbf{X}(\mathcal{T})$ in the pointwise sense. Hence $ S^\wedge(c)\in\mathcal{C}.$\\

In this way we associate with each operator
  \(S\in\mathbf{H}(\mathcal{T})\)  a mapping $c\mapsto S^\wedge(c)$ from $\mathbb{B}(\mathcal{T})$ to the algebra of
  constants $\mathcal{C}$. We shall refer to this mapping as {\em
  the  point evaluation of $S$.} Its main properties are listed in the
 following lemma.

\begin{lemma}\label{poevprop}\ \\
\begin{enumerate}
\item
  Let $F,G\in\mathbf{H}(\mathcal{T})$, ${p}\in{\mathcal{C}}$, ${c}\in\mathbb{B}(\mathcal{T})$.
   Then
\begin{align}
\label{linpoiev}(F{p}+G)^{\wedge}({c})&=F^{\wedge}
    ({c})\cdot{p}+G^{\wedge}({c}),\\
 \label{mulpoiev} (FG)^{\wedge}({c})&=(F^{\wedge}({c})\cdot G)^{\wedge}({c}).
  \end{align}

  \item If $S\in\mathbf{H}(\mathcal{T})$ and  $S^\wedge
({c})=0$ for every ${c}\in\mathbb{B}(\mathcal{T})$, then $S=0.$
\end{enumerate}
  \end{lemma}
  \begin{proof}\
\paragraph{(I)} In view of Remark \ref{linexp1}, the identity \eqref{linpoiev} follows immediately from the definition
\eqref{lepoev} of the point evaluation. Therefore, it suffices to establish the identity \eqref{mulpoiev} for
$F,G$ of the form
\[F=w_1^*,\ G=w_2^*,\quad\text{where } w_1,w_2\in\mathcal{F}_q.\]
But Lemma \ref{trivlem2} implies that
\begin{multline*}(w_1^*w_2^*)^\wedge
({c})=(c\alpha)^{(|w_1|+|w_2|)}w_1^*w_2^*=
(c\alpha)^{|w_2|}(w_1^*)^\wedge(c)w_2^*\\=\sum_{|w|=|w_2|}(c\alpha)^{|w|}w^*
\left(w(w_1^*)^\wedge(c)w_2^*\right)
=\sum_{|w|=|w_2|}(w^*)^\wedge(c)\left(w(w_1^*)^\wedge(c)w_2^*\right)\\
=\left(\sum_{|w|=|w_2|}w^*w(w_1^*)^\wedge(c)w_2^*\right){}^\wedge({c})= ((w_1^*)^\wedge(c)w_2^*)^\wedge({c}).
\end{multline*}
\ \\

\paragraph{(II)} We prove that $S_{[w]}=0$ for each $w\in\mathcal{F}_q$. We use induction on $|w|$. First,
\[S_{[\oo]}=S^\wedge(0)=0.\] Next, assume that
\[S_{[w]}=0\quad \forall w\,:\,|w|\leq n\]
and let \[w_{n+1}=\alpha_{i_1}\alpha_{i_2}\cdots\alpha_{i_{n+1}}.\] We shall prove that $S_{[w_{n+1}]}=0.$\\

Denote
$$w_k=\alpha_{i_1}\dots\alpha_{i_k},\quad 1\leq k\leq n+1,\quad w_0=\oo.$$
 Fix $t_0\in\mathcal{T}$ and
 consider the diagonal operators $c_j\in\mathcal{C},$ $1\leq j\leq q,$ defined as follows:
\[c_j\chi_t=\left\{\begin{aligned}\chi_{t_0w_k},&\quad\text{if}\quad  j=i_{k+1}\text{ and
}t=t_0w_k,\quad 0\leq k\leq n,\\
0,&\quad\text{otherwise}.\end{aligned} \right.\] We set
\[c=\begin{pmatrix}c_1&\dots &c_q\end{pmatrix}\]
and observe that, as follows from \eqref{prsh},
\[
(c\alpha)\chi_t=\left\{\begin{aligned}\chi_{t_0w_k},&\quad\text{if}
\quad t=t_0w_{k+1},\quad 0\leq k\leq n,\\
0,&\quad\text{otherwise}.\end{aligned}\right.\] Hence
\begin{align*}
(c\alpha)^{n+1}\chi_t&=\left\{\begin{aligned}\chi_{t_0},&\quad\text{if}
\quad t=t_0w_{n+1},\\ 0,&\quad\text{otherwise},\end{aligned}\right.\\
(c\alpha)^{n+2}&=0.
\end{align*}
In particular,  \(c\in\mathbb{B}(\mathcal{T}).\)\\

 Furthermore, in view of \eqref{fqw}, for
$w\in\mathcal{F}_q$ such that $|w|=n+1$ we have
\[(c\alpha)^{n+1}w^*\chi_t=(c\alpha)^{n+1}\chi_{tw}=\left\{\begin{aligned}\chi_{t_0}&\text{ if
}t=t_0\text{ and }w=w_{n+1},\\
0,&\text{ otherwise}.\end{aligned}\right.\] Thus by the induction assumption
\[0=S^\wedge(c)=(c\alpha)^{n+1}w_{n+1}^*S_{[w_{n+1}]},\]
which implies (see \eqref{fqw})
\[0=S_{[w_{n+1}]}^*w_{n+1}(\alpha^*c^*)^{n+1}\chi_{t_0}
=S_{[w_{n+1}]}^*w_{n+1}\chi_{t_0w_{n+1}}=S_{[w_{n+1}]}^*\chi_{t_0}.\] Since
 $t_0\in\mathcal{T}$ was chosen arbitrarily, $S_{[w_{n+1}]}=0.$
  \end{proof}

\section{The  space of causal Hilbert -- Schmidt operators}\label{s5}
In this section we consider the following spaces of Hilbert -- Schmidt operators:
\begin{align*}\mathbf{H}_2(\mathcal T)&\defi\{F\in\mathbf{H}(\mathcal{T});\,\|F\|_2^2\defi
\trace( F^*F)<\infty\},\\
\mathcal{C}_2&\defi\mathcal{C}\cap\mathbf{H}_2(\mathcal T).\end{align*}

As we shall see from Propositions \ref{cla1} and \ref{hilpoex} below,  the space $\mathbf{H}_2(\mathcal T)$ of
causal Hilbert -- Schmidt operators plays the role of the Hardy space of the unit disk in the present setting:
elements of the algebra $\mathbf{H}(\mathcal T)$ act on the space $\mathbf{H}_2(\mathcal T)$ by multiplication.
The space $\mathcal{C}_2$ is the space of constants which appear in the power series expansions (see Theorem
\ref{trivprop}) of the elements of $\mathbf{H}_2(\mathcal T)$.

\begin{proposition}\label{cla1}\ \\
\begin{enumerate}\item
The space $\mathbf{H}_2(\mathcal{T})$ is a Hilbert space
 contractively included in
\(\mathbf{H}(\mathcal T)\):
\[\forall
F\in\mathbf{H}_2(\mathcal{T})\quad \langle F,F\rangle_2\leq \|F\|^2.\]

\item  $\mathcal{C}_2$ is a closed subspace of the Hilbert space $\mathbf{H}_2(\mathcal{T})$.

\item For every \(F\in\mathbf{H}_2(\mathcal{T})\) and \(S\in\mathbf{H}(\mathcal{T})\) the operators
 \(SF\) and \(FS\) belong to \(\mathbf{H}_2(\mathcal{T}) \)
and it holds that \begin{equation}\label{trivest2}\max(\|SF\|_2,\|FS\|_2)\leq\|S\|\,\|F\|_2.\end{equation}
Moreover, the multiplication operators $M_S,\hat{M}_S$ defined on \(\mathbf{H}_2(\mathcal{T})\) by
\begin{equation}\label{defmul}
M_SF\defi SF,\quad \hat{M}_SF\defi FS, \quad F\in\mathbf{H}_2(\mathcal{T}),\end{equation} satisfy
\[\|M_S\|=\|\hat{M}_S\|=\|S\|.\]
\end{enumerate}
\end{proposition}

\begin{proof}
It is  well known (see, for instance, \cite{gk-1}) that the space of Hilbert -- Schmidt operators on a given
separable Hilbert space is a Hilbert space. In particular, the space \[\mathbf{X}_2(\mathcal
T)\defi\{F\in\mathbf{X}(\mathcal{T});\,\|F\|_2^2\defi \trace( F^*F)<\infty\}\]  is a Hilbert space. For every
$F\in\mathbf{X}_2(\mathcal{T})$ and $f\in\ell_2(\mathcal{T})$ it holds that
\begin{multline*}
\|Ff\|_{\ell_2(\mathcal{T})}^2=\sum_{t\in\mathcal{T}}|Ff(t)|^2=\sum_{t\in\mathcal{T}}\left|\sum_{s\in\mathcal{T}}
F\chi_s(t)\cdot f(s)\right|^2\\\leq\sum_{t,s\in\mathcal{T}}|
F\chi_s(t)|^2\,\|f\|_{\ell_2(\mathcal{T})}^2=\|F\|_2^2\, \|f\|_{\ell_2(\mathcal{T})}^2,\end{multline*} hence the
space \(\mathbf{X}_2(\mathcal T)\) is contractively included in \(\mathbf{X}(\mathcal T).\)\\

The space $\mathbf{H}_2(\mathcal{T})$ is the intersection
\begin{equation}\label{cauint}
\mathbf{H}_2(\mathcal{T})=\mathbf{X}_2(\mathcal T)\cap\mathbf{H}(\mathcal{T}).\end{equation} Since the algebra
$\mathbf{H}(\mathcal{T})$ is closed in $\mathbf{X}(\mathcal{T})$ in the pointwise sense, it is also closed with
respect to the operator norm. It follows that  $\mathbf{H}_2(\mathcal{T})$ is a closed (with respect to the
Hilbert -- Schmidt norm $\|\cdot\|_2$) subspace of the Hilbert space $\mathbf{X}_2(\mathcal T)$,
which proves the statement (I).\\

The proof of the statement (II) is analogous: it uses the fact that the algebra of constants $\mathcal{C}$ is
closed
in $\mathbf{H}(\mathcal{T})$ in the pointwise sense and hence also  with respect to the operator norm.\\

In order to prove the statement (III), we observe that for every $S\in\mathbf{X}(\mathcal{T})$ and
$F\in\mathbf{X}_2(\mathcal{T})$ it holds that
\begin{align*}\|SF\|_2^2&=\sum_{t\in\mathcal{T}}\|SF\chi_t\|_{\ell_2(\mathcal{T})}^2 \leq
\|S\|^2\sum_{t\in\mathcal{T}}\|F\chi_t\|_{\ell_2(\mathcal{T})}^2=\|S\|^2\,\|F\|_2^2;\\
\|FS\|_2^2&=\|S^*F^*\|_2^2\leq\|S^*\|^2\,\|F^*\|_2^2=\|S\|^2\,\|F\|_2^2.\end{align*}

In particular, taking into account \eqref{cauint} and the fact that $\mathbf{H}(\mathcal{T})$ is an algebra, we
may conclude that \[\forall S\in\mathbf{H}(\mathcal{T}), \forall F\in\mathbf{H}_2(\mathcal{T})\quad
SF,FS\in\mathbf{H}_2(\mathcal{T})\] and \eqref{trivest2} holds. \\

Let $f\in\spa\{\chi_t\}$ and choose  $t_0\in\mathcal{T}$  such that $f(t)=0$ $\forall t\preceq t_0.$ Consider
the operator $F$ defined by \begin{equation}\label{defff}Fu= u(t_0)\cdot f,\quad
u\in\ell_2(\mathcal{T}).\end{equation} Then
\[\|F\|_2^2=\sum_{t\in\mathcal{T}}\|F\chi_t\|_{\ell_2(\mathcal{T})}^2=\|f\|_{\ell_2(\mathcal{T})}^2,\]
hence $F\in\mathbf{X}_2(\mathcal{T}).$  According to Definition \eqref{acau}, the operator $F$ is causal, hence
$F\in\mathbf{H}_2(\mathcal{T}).$\\

 Furthermore, let $S\in\mathbf{H}(\mathcal{T}).$ Then
\[\|SF\|_2^2=\sum_{t\in\mathcal{T}}\|SF\chi_t\|_{\ell_2(\mathcal{T})}^2=\|Sf\|_{\ell_2(\mathcal{T})}^2,\]
Since $\spa\{\chi_t\}$ is dense in $\ell_2(\mathcal{T}),$ we conclude that the left multiplication operator
$M_S$ satisfies $\|M_S\|\geq\|S\|.$ On the other hand, \eqref{trivest2} implies $\|M_S\|\leq\|S\|,$ hence
$\|{M}_S\|=\|S\|.$\\

 The proof of the equality $\|\hat{M}_S\|=\|S\|$ for the right multiplication operator
$\hat{M}_S$ is analogous.
\end{proof}

\begin{proposition}\label{hilpoex}
 Let $F\in\mathbf{H}(\mathcal{T})$ and let
\begin{equation}\label{llg5}
F=\sum_{w\in\mathcal{F}_q}w^*F_{[w]},\end{equation} where $F_{[w]}\in\mathcal{C},$ be the pointwise absolutely
convergent expansion for $F$, as in Theorem \ref{trivprop}. Then $F\in\mathbf{H}_2(\mathcal T)$ if and only if
\begin{equation}\label{hardy}\forall w\in\mathcal{F}_q\
F_{[w]}\in\mathcal{C}_2\quad\text{and}\quad\sum_{w\in\mathcal{F}_q}\|F_{[w]}\|_2^2<\infty.\end{equation}
 In
this case the expansion \eqref{llg5} converges in the $\mathbf{H}_2(\mathcal{T})$-norm and
\begin{equation}\label{hardynorm}\|F\|_2^2=\sum_{w\in\mathcal{F}_q}\|F_{[w]}\|_2^2.\end{equation}
\end{proposition}

\begin{proof}\ \paragraph{$(\Rightarrow)$}
Assume  that $F\in\mathbf{H}_2(\mathcal T)$. Then, as follows from the relations \eqref{tayco} in Theorem
\ref{trivprop} and the statement (III) of Proposition \ref{cla1}, we have
\begin{multline*}\|F\|_2^2=\|w\|\,\|F\|_2^2\geq\|wF\|_2^2 =\sum_{t,s\in\mathcal{T}}| (wF\chi_t)(s)|^2\\
\geq \sum_{\substack{t,s\in\mathcal{T}\\t\asymp s}}| (wF\chi_t)(s)|^2= \sum_{t,s\in\mathcal{T}}|
(F_{[w]}\chi_t)(s)|^2=\|F{[w]}\|_2^2,\end{multline*} hence
$F_{[w]}\in\mathcal{C}_2$.\\

Furthermore, since $F$ is causal,
\begin{multline}\label{trivcomp1}
\|F\|_2^2=\sum_{\substack{s,t\in\mathcal{T}\\s\preceq t}}|
(F\chi_s)(t)|^2=\sum_{\substack{s,t\in\mathcal{T}\\s\asymp
t}}\sum_{w\in\mathcal{F}_q}|(F\chi_s)(tw)|^2\\=\sum_{\substack{s,t\in\mathcal{T}\\s\asymp
t}}\sum_{w\in\mathcal{F}_q}|(F_{[w]}\chi_s(t)|^2=\sum_{w\in\mathcal{F}_q}\|F_{[w]}\|_2^2.
\end{multline}
Thus \eqref{hardynorm} holds.\\

 As a consequence, we obtain the convergence of the expansion \eqref{llg5} in
the following sense: if we order somehow the countable set $\mathcal{F}_q$, say
$$\mathcal{F}_q=\{w_j\}_{j=0}^\infty,$$
then, according to \eqref{hardy} and the statement (III) of Proposition \ref{cla1}, for each $n=0,1,2,\dots$ the
finite sum $\sum_{j=0}^nw_j^*F_{[w_j]}$ belongs to $\mathbf{H}_2(\mathcal{T})$ and it holds that
\[\left\|F-\sum_{j=0}^nw_j^*F_{[w_j]}\right\|_2^2=\sum_{j=n+1}^\infty\|F_{[w_j]}\|_2^2\rightarrow 0\text{
as }n\rightarrow\infty.
\]

\paragraph{$(\Leftarrow)$} Assume that $F\in\mathbf{H}(\mathcal{T})$ and
that the coefficients $F_{[w]}$ of the expansion
\eqref{llg5} satisfy the conditions \eqref{hardy}. Then it suffices to reverse
 the computation \eqref{trivcomp1}
in order to show that $F\in\mathbf{H}_2(\mathcal{T}).$
\end{proof}

Our next goal is to demonstrate that the space $\mathbf{H}_2(\mathcal{T})$ has a reproducing kernel structure
with respect to the point evaluation defined in the previous section (see \eqref{lepoev}).

\begin{lemma}\label{cor} Let $F\in\mathbf{H}_2(\mathcal{T})$ have the expansion
\eqref{llg5}, where $F_{[w]}\in\mathcal{C}_2,$ and let $c\in\mathbb{B}(\mathcal{T}).$ Then the series
\begin{equation}\label{lpoe}
F^\wedge(c)=\sum_{n=0}^\infty(c\alpha)^n\left(\sum_{\substack{w\in\mathcal{F}_q\\|w|=n}}w^*F_{[w]}\right)\end{equation}
converges absolutely in $\mathcal{C}_2.$
\end{lemma}
\begin{proof}
It follows from the identity \eqref{hardynorm} in Proposition \ref{hilpoex} that
\[\left\|\sum_{\substack{w\in\mathcal{F}_q\\|w|=n}}w^*F_{[w]}\right\|_2\leq\|F\|_2,\quad n=0,1,2\dots.\]
 Hence, in view of the definition \eqref{ball} of $\mathbb{B}(\mathcal{T})$ and the estimate \eqref{trivest2}
in Proposition \ref{cla1}, the series  \eqref{lpoe} converges absolutely in $\mathbf{H}_2(\mathcal{T})$. Since
each term of the series belongs to $\mathcal{C}$, the desired conclusion follows.
\end{proof}

\begin{theorem}\label{proker}
Let $c\in\mathbb{B}(\mathcal{T})$. Then the operator $I-{\alpha^*} {c}^*$ is invertible in
$\mathbf{H}(\mathcal{T})$ and its inverse
\begin{equation}
\label{rk}K^{{c}}_{\wedge}\defi (I-{\alpha^*} {c}^*)^{-1}
\end{equation}
 satisfies
\begin{equation}\label{al}
    \langle F^{\wedge}({c}),{ k}\rangle_2=\langle F,
    K_{\wedge}^{ c}\,{ k}\rangle_2,\quad \forall
    F\in\mathbf{H}_2(\mathcal{T}),\forall k\in\mathcal{C}_2.
   \end{equation}
\end{theorem}
\begin{proof}
In view of the definition \eqref{ball} of $\mathbb{B}(\mathcal{T}),$ $K^c_\wedge$ is the sum of the absolutely
convergent series
 \begin{equation}\label{kerser}K^c_\wedge=\sum_{n=0}^\infty(\alpha^* c^*)^n.\end{equation}
Since each term of the series belongs to $\mathbf{H}(\mathcal{T})$, so does $K^c_\wedge$. Let us now choose and
fix  an element $k\in\mathcal{C}_2$ and an element $F\in\mathbf{H}_2(\mathcal{T})$ with the expansion
\eqref{llg5}, where $F_{[w]}\in\mathcal{C}_2.$
 Then, according to the statement (III) of Proposition
\ref{cla1},
$$K^c_\wedge\,k=\sum_{n=0}^\infty(\alpha^* c^*)^nk\in\mathbf{H}_2(\mathcal{T}),$$ where  convergence is absolute with
respect to the $\mathbf{H}_2(\mathcal{T})$-norm. Since, in view of the identity  \eqref{trivid},
\begin{equation}\label{acn}(\alpha^*c^*)^n=\sum_{\substack{w\in\mathcal{F}_q\\|w|=n}}w^*w(\alpha^*
c^*)^n=\sum_{\substack{w\in\mathcal{F}_q\\|w|=n}}w^*(w^*)^\wedge(c)^*,\end{equation} Proposition \ref{hilpoex}
and Lemma \ref{cor} imply that
\begin{multline*}
\langle F,K^c_\wedge\,k\rangle_2=\sum_{n=0}^\infty\langle F,(\alpha^*c^*)^nk\rangle_2=
\sum_{n=0}^\infty\langle\sum_{\substack{w\in\mathcal{F}_q\\|w|=n}}
w^*F_{[w]},(\alpha^*c^*)^nk\rangle_2\\=\sum_{n=0}^\infty\langle(c\alpha)^n\sum_{\substack{w\in\mathcal{F}_q\\|w|=n}}
w^*F_{[w]},k\rangle_2=\langle F^\wedge(c),k\rangle_2.
\end{multline*}
\end{proof}
\begin{remark}\label{dense} Note that Theorem
\ref{proker} and the statement (II) of Lemma \ref{poevprop}   imply that
$$\overline{\spa}\{K^c_\wedge\,k\,:\,c\in\mathbb{B}(\mathcal{T}),k\in\mathcal{C}_2\}=\mathbf{H}_2(\mathcal{T}).$$
\end{remark}

We close this section with the description of the counterparts of the backward shift operator $R_0$ in the
stationary single-scale setting (see the formula \eqref{defr0} in Introduction).

\begin{proposition}\label{glele} Let operators
$A_j:\mathbf{H}_2(\mathcal{T})\longrightarrow\mathbf{H}_2(\mathcal{T}),$ $j=1,\dots, q,$ be defined by
\begin{equation}\label{gl2}
A_jF\defi(F-F^\wedge(0))\alpha_j,\quad F\in\mathbf{H}_2(\mathcal{T}),\ j=1,\dots, q.\end{equation}
 Then:
 \begin{enumerate}
\item The operator $A_j$ is
  the adjoint of the right multiplication operator
\(\hat{M}_{\alpha_j^*}\) in $\mathbf{H}_2(\mathcal{T}):$
\begin{equation*}
A_j=\hat{M}_{\alpha_j^*}^*,\quad j=1,\dots, q.\end{equation*}

\item For $i,j=1,\dots, q$ the following relations hold:
\begin{equation*}
 A_j\hat{M}_{\alpha_i^*}=\hat{M}_{\alpha_i^*\alpha_j},
\qquad\hat{M}_{\alpha_i^*}A_j=\left\{\begin{aligned}(I-C^*C),&\quad\text{if}\quad i=j,\\
0,&\quad\text{otherwise,}\end{aligned}\right.
\end{equation*}
where
\begin{equation}\label{ev0}
CF\defi F^\wedge(0).\end{equation}

\item Let $F\in\mathbf{H}_2(\mathcal{T})$ with the expansion \eqref{llg5}, where $F_{[w]}\in\mathcal{C}_2,$ be
given. Then for every pair $w,v\in\mathcal{F}_q$ such that $|v|=|w|$ it holds that
\begin{equation*}
F_{[w]}=w\cdot(CA^vF)\cdot v^*,
\end{equation*}
 where
\begin{equation}\label{stnot3}
A^{\alpha_{i_1}\dots\alpha_{i_k}}\defi A_{i_k}\dots A_{i_1},\quad A^{\oo}\defi I.
\end{equation}
\end{enumerate}
\end{proposition}
\begin{proof}
Consider $F\in\mathbf{H}_2(\mathcal{T})$ with the expansion
\eqref{llg5}, where $F_{[w]}\in\mathcal{C}_2.$
 Then, in view of  Proposition \ref{hilpoex},
$$A_jF=(F-F_{[\oo]})\alpha_j=\sum_{\substack{w\in\mathcal{F}_q\\|w|\geq 1}}w^*F_{[w]}\alpha_j
\in\mathbf{H}_2(\mathcal{T})$$ and for every $G\in\mathbf{H}_2(\mathcal{T})$ it holds that
\begin{multline*}\langle F,G\alpha_j^*\rangle_2=\langle F-F_{[\oo]}, G\alpha_j^*\rangle_2=\trace(\alpha_j
G^*(F-F_{[\oo]}))\\=\trace( G^*(F-F_{[\oo]})\alpha_j)=\langle A_jF,G\rangle_2.\end{multline*} This proves the
statement (I) of the Proposition.\\

 The statement (II) follows immediately from the Cuntz relations
\eqref{cuntz}.\\

 In order to prove the statement (III), let us fix a sequence of indices $i_1,i_2,i_3,\dots,$ $1\leq i_n\leq q.$  Then
\begin{multline*}F=F_{[\oo]}+(F-F_{[\oo]})=CF+(A_{i_1}F)\alpha_{i_1}^*\\=CF+(CA_{i_1}F)\alpha_{i_1}^*+(A_{i_2}A_{i_1}F)
\alpha_{i_2}^*\alpha_{i_1}^*=\dots\\ =\sum_{n=0}^m (CA^{v_n}F)v_n^*+(A^{v_{m+1}}F)v_{m+1}^*,\quad
 m=0,1,2,\dots\end{multline*}
 where $v_0=\oo,$
$v_n=\alpha_{i_1}\dots\alpha_{i_n}.$ Now it follows from the identity \eqref{tayco} in
Theorem \ref{trivprop}
that
$$\sum_{|w|=n}w^*F_{[w]}d=(CA^{v_n}F)v_n^*\quad \forall n$$
and hence
$$F_{[w]}=w(CA^{v_{|w|}}F)v_{|w|}^*\quad \forall w.$$
Since the sequence of indices $i_1,i_2,i_3,\dots$ was chosen arbitrarily,
this completes the proof.
\end{proof}

\section{Schur multipliers and de Branges -- Rovnyak spaces}\label{s6}
\begin{definition}
Let \(S\in\mathbf{H}(\mathcal{T})\) be such that \(\|S\|\leq 1.\) Then \(S\)  is said to be a {\em Schur
multiplier.}
\end{definition}

\begin{theorem} \label{Schthrm}
Let a mapping $s:\mathbb{B}(\mathcal{T})\mapsto\mathcal{C}$ be given. Then there exists a Schur multiplier
\(S\in\mathbf{H}(\mathcal{T})\) such that
\[s(c)=S^\wedge(c)\quad \forall c\in\mathbb{B}(\mathcal{T})\]
   if and only if the kernel
$K_s:{\mathbb{B}}({\mathcal{T}})\times{\mathbb{B}}({\mathcal{T}})\to{\mathcal{C}}$ defined by
  \begin{equation}\label{defks} K_s( c, d)\defi
\sum_{n=0}^\infty (c\alpha)^n(I-s(c)s(d)^*)(d\alpha)^{n*},\quad c,d\in\mathbb{B}(\mathcal{T}),
\end{equation}
  is positive: for any $m\geq 0$, $c_0,...,c_m
\in\mathbb{B}(\mathcal{T})$, $k_0,\dots,k_m\in\mathcal{C}_2$,
  it holds that
  \[ \sum_{i,j=0}^{m}\langle K_s(c_i,
c_j)k_j,k_i\rangle_2\geq 0.\]

In this case  $K_s(c,d)= (K_S^d)^\wedge(c),$ where
\begin{equation}\label{repkers}K_S^{d}\defi(I-SS^\wedge({d})^*)K^{d}_\wedge,
\quad d\in\mathbb{B}(\mathcal{T}).\end{equation}
  \end{theorem}
In the proof of Theorem \ref{Schthrm} we shall use the following lemma:
\begin{lemma}\label{commut}\ \\
\begin{enumerate}
\item
Let $T:\mathbf{H}_2(\mathcal{T})\longrightarrow\mathbf{H}_2(\mathcal{T})$ be a bounded linear operator.
 Then $T$
is of the form $T=M_S$ for some $S\in\mathbf{H}(\mathcal{T})$ if and only if
\[\forall Q\in\mathbf{H}(\mathcal{T})\quad T\hat{M}_Q=\hat{M}_QT.\]

\item Let $P:\mathcal{C}_2\longrightarrow\mathcal{C}_2$ be a bounded linear operator.
 Then $T$
is of the form $T=M_c$ for some $c\in\mathcal{C}$ if and only if
\[\forall d\in\mathcal{C}\quad P\hat{M}_d=\hat{M}_dP.\]
\end{enumerate}
\end{lemma}
\begin{proof}
We shall prove only the statement (I) of the Proposition; the proof of the statement (II) is completely
analogous.\\

\paragraph{$(\Rightarrow)$} The "only if" direction is clear: for every $S,Q\in\mathbf{H}(\mathcal{T})$ and
$F\in\mathbf{H}_2(\mathcal{T})$ it holds that
\[M_S\hat{M}_QF=SFQ=\hat{M}_QM_SF.\]
\ \\

\paragraph{$(\Leftarrow)$} Assume that an operator $T$, which commutes with every $\hat{M}_Q$, is given.
For each $t\in\mathcal{T}$ let us consider the projection $\pi_t\in\mathcal{C}_2$,
$$\pi_tu\defi u(t)\chi_t,\quad u\in\ell_2(\mathcal{T}),$$ and define an operator $S$ on
$\spa\{\chi_t\}$ by
$$S\chi_t=(T\pi_t)\chi_t.$$

Let $f\in\spa\{\chi_t\}$ and choose  $t_0\in\mathcal{T}$  such that $f(t)=0$ $\forall t\preceq t_0.$ Consider
$F\in\mathbf{H}_2(\mathcal{T})$ defined by \[Fu= u(t_0)\cdot f,\quad u\in\ell_2(\mathcal{T}),\] as in the proof
of Proposition \ref{cla1} (see \eqref{defff}). Then
$$F=\sum_{t\in\mathcal{T}}f(t)\cdot F_t^{t_0},\quad\text{where}\quad F_t^{t_0}u=
u(t_0)\chi_t,\quad u\in\ell_2(\mathcal{T}).$$ Note that $F_t^{t_0}=\pi_t
F_t^{t_0}=F_t^{t_0}\pi_{t_0}\in\mathbf{H}_2(\mathcal{T})$, hence for every $u\in\ell_2(\mathcal{T})$ we have
\begin{multline*}(TF)u=u(t_0)\sum_{t\in\mathcal{T}}f(t)(TF_t^{t_0})\chi_{t_0}
=u(t_0)\sum_{t\in\mathcal{T}}f(t)(T\pi_t)F_t^{t_0}\chi_{t_0}\\
=u(t_0)\sum_{t\in\mathcal{T}}f(t)(T\pi_t)\chi_{t}=u(t_0)Sf=SFu.\end{multline*} It follows that
$$\|Sf\|_{\ell_2(\mathcal{T})}\leq \|TF\|\leq \|TF\|_2\leq\|T\|\cdot\|F\|_2=\|T\|\cdot
\|f\|_{\ell_2(\mathcal{T})}$$  and $(Sf)(t)=0$ $\forall t\preceq t_0.$ Since
$$\overline{\spa}\{\chi_t\,:\,t\in\mathcal{T}\}=\ell_2(\mathcal{T})
\quad\text{and}\quad\overline{\spa}\{F_t^s\,:\,t,s\in\mathcal{T},s\preceq t\}
 =\mathbf{H}_2(\mathcal{T}),$$ we conclude that
 $$S\in\mathbf{H}(\mathcal{T})\quad\text{and}\quad T=M_S.$$
\end{proof}

\begin{proof}[Proof of Theorem \ref{Schthrm}]
Let us assume first that $s(c)=S^\wedge(c),$ where $S$ is a Schur multiplier. According to  Theorem
\ref{proker}, for $k\in\mathbf{C}_2$ and $F\in\mathbf{H}_2(\mathcal{T})$ we have
\[\langle K^d_\wedge\,k,SF\rangle_2 =\langle
K^d_\wedge\,k,S^\wedge(d)F\rangle_2=\langle S^\wedge(d)^*K^d_\wedge\,k,F\rangle_2,\] hence
\begin{align}\label{msta}
M_S^*(K^d_\wedge\,k)&=S^\wedge(d)^*K^d_\wedge\,k,\\
\label{mstb}K_S^d\,k&=(I-M_SM_S^*)(K^d_\wedge\,k).
\end{align}
Furthermore, it follows from the identities \eqref{kerser} and \eqref{acn}, the definition \eqref{lepoev} of the
point evaluation and the statement (I) in Lemma \ref{poevprop} that
\begin{multline*}(K^d_S)^\wedge(c)=
(K^d_\wedge)^\wedge(c)-(s(c)s(d)^*K^d_\wedge)^\wedge(c)\\=
\sum_{n=0}^\infty(c\alpha)^n(I-s(c)s(d)^*)(\alpha^*d^*)^n =K_s(c,d).
\end{multline*}

 Now, given $c_0,...,c_m
\in\mathbb{B}(\mathcal{T})$, $k_0,\dots,k_m\in\mathcal{C}_2$, we observe that  \[ \sum_{i,j=0}^{m}\langle
K_s(c_i, c_j)k_j,k_i\rangle_2=\left\langle
(I-M_SM_S^*)\left(\sum_{j=0}^mK_\wedge^{c_j}k_j\right),\sum_{i=0}^mK_\wedge^{c_i}k_i\right\rangle_2\geq 0,\]
because  the operator $I-M_SM_S^*$ is
positive. Thus the kernel $K_s(c,d)$ is positive.\\

 Conversely, assume that the kernel $K_s(c,d)$ is positive and define
on $\spa\{K^d_\wedge k\,:\,d\in\mathbb{B}(\mathcal{T}),k\in\mathcal{C}_2\}$ an operator $T$ by
$$T(K^d_\wedge
k)=s(d)^*K^d_\wedge k.$$ Then $T$ is a well-defined contraction, because
\[\left\|T\sum_{j=0}^m K^{c_j}_\wedge\,k_j\right\|_2^2=\left\|\sum_{j=0}^m K^{c_j}_\wedge\,k_j\right\|_2^2-
\sum_{i,j=0}^{m}\langle K_s(c_i, c_j)k_j,k_i\rangle_2\leq\left\|\sum_{j=0}^m K^{c_j}_\wedge\,k_j\right\|_2^2.\]
Hence, in view of Remark \ref{dense}, $T$ can be extended as a contraction on $\mathbf{H}_2(\mathcal{T})).$
\\

The adjoint operator has the property \begin{equation}\label{tsta}(T^*F)^\wedge(c)=(s(c)F)^\wedge(c)\quad
\forall F\in\mathbf{H}_2(\mathcal{T}), \forall c\in\mathbb{B}(\mathcal{T}).\end{equation} In view of the
statement (I) in Lemma \ref{poevprop}, for every $Q\in\mathbf{H}(\mathcal{T}),$ $F\in\mathbf{H}_2(\mathcal{T}),$
$c\in\mathbb{B}(\mathcal{T})$ we have
\begin{multline*}(T^*\hat{M}_QF)^\wedge(c)=(s(c)FQ)^\wedge(c)=((s(c)F)^\wedge(c)Q)^\wedge(c)\\
=((T^*F)^\wedge(c)Q)^\wedge(c)=(\hat{M}_QT^*F)^\wedge(c),\end{multline*} which, according to the statement (II)
of the same Lemma \ref{poevprop}, implies
\[T^*\hat{M}_Q=\hat{M}_QT^* \forall Q\in\mathbf{H}(\mathcal{T}).\]
  Now it follows from Lemma \ref{commut} that there exists a Schur multiplier $S$
such that $T^*=M_S.$ In view of \eqref{tsta} and Lemma \ref{poevprop}, this Schur multiplier $S$ satisfies
 $$S^\wedge(c)=s(c)\quad\forall
c\in\mathbb{B}(\mathcal{T}).$$
\end{proof}

We recognize in kernel \eqref{defks} the analogue of the kernel \eqref{00ker}, mentioned in Introduction. As in
the single-scale case, we consider the associated de Branges -- Rovnyak space defined below.

\begin{definition}\label{debro} Let $S\in\mathbf{H}(\mathcal{T})$ be a Schur multiplier, let
\[\mathcal{B}_S\defi I-M_SM_S^*\]
and let \(\pi_S\) denote the orthogonal projection in \({\bf H}_2({\mathcal T})\) onto \(\ker\mathcal{B}_S.\)
 The
Hilbert space \(\mathbf{H}(S)\), defined by
\begin{equation*}
\mathbf{H}(S)=\sqrt{\mathcal{B}_S}\,\mathbf{H}_2(\mathcal{T});\quad \|
\sqrt{\mathcal{B}_S}\,F\|_{\mathbf{H}(S)}= \|(I-\pi_S)F\|_{\mathbf{H}_2(\mathcal{T})},
\end{equation*}
 is said to be the {\em de Branges -- Rovnyak space} associated with
\(S\).
\end{definition}
In the sequel we shall use the following terminology:
\begin{definition}\label{rightc}
Let $\mathbf{H}$ be a Hilbert space of elements of $\mathbf{H}_2(\mathcal{T})$. The space $\mathbf{H}$ is said
to be {\em right $\mathcal{C}$-invariant} if for every $F,G\in\mathbf{H}$ and $c\in\mathcal{C}$ it holds that
\[Fc\in\mathbf{H},\quad\|Fc\|_{\mathbf{H}}\leq\|F\|_{\mathbf{H}}\,\|c\|,
\quad\langle Fc,G\rangle_{\mathbf{H}}=\langle F,Gc^*\rangle_{\mathbf{H}}.\]
\end{definition}
\begin{proposition}\label{cla2}
Let $S\in\mathbf{H}(\mathcal{T})$ be a Schur multiplier. Then  the de Branges -- Rovnyak space
\(\mathbf{H}(S)\), associated with $S$, is a right $\mathcal{C}$-invariant Hilbert space.

 Furthermore,
for every $F\in\mathbf{H}(S),$ ${k}\in\mathcal{C}_2,$ $c\in\mathbb{B}(\mathcal{T})$ it holds that
\begin{equation} \label{repdbr} K_S^{c}{k} \in\mathbf{H}(S)\quad\text{and}\quad\langle F,
K_S^{c}{k}\rangle_{\mathbf{H}(S)}=\langle F^\wedge ({c}),{k}\rangle_2,
\end{equation}
where $K_S^{c}$ is as in \eqref{repkers}. In particular,
\begin{equation}\label{dense2}
\mathbf{H}(S)=\overline{\spa}\{K_S^{c}{k}\,:\,{k}\in\mathcal{C}_2,c\in\mathbb{B}(\mathcal{T})\}.
\end{equation}
\end{proposition}
\begin{proof}
Let $c\in\mathcal{C}$. Since  the adjoint of the right multiplication operator $\hat{M}_c$ in
$\mathbf{H}_2(\mathcal{T})$ is given by $\hat{M}_c^*=\hat{M}_{c^*},$ Lemma \ref{commut} implies that the
operators $\hat{M}_c$ and $\mathcal{B}_S$ commute. Hence, in view of Definitions \ref{debro}, \ref{rightc} and
the statement (III) of Proposition \ref{cla1}, the space $\mathbf{H}(S)$ is right $\mathcal{C}$-invariant. \\

Furthermore, as follows from  \eqref{mstb} and Definition \ref{debro}, for every $k\in\mathcal{C}_2$
$K_S^ck\in\mathbf{H}(S)$. Moreover, for every $F\in\mathbf{H}(S)$ we have
\[\langle F,K_S^ck\rangle_{\mathbf{H}(S)}=\langle F, K^c_\wedge\,k\rangle_2.\]
Thus the identity \eqref{al} in Theorem \ref{proker} implies \eqref{repdbr} and the statement (II) in Lemma
\ref{poevprop} implies \eqref{dense2}.
\end{proof}

\begin{theorem}\label{main1}
Let \(S\in\mathbf{H}(\mathcal{T})\) be a Schur multiplier and let $\mathbf{H}(S)$ be the associated de Branges
-- Rovnyak space. Set
\begin{equation}\label{defabcd}\begin{split} A_jF=(F-F^\wedge(0))\alpha_j,&\quad B_jd=(S-S^\wedge(0))d\alpha_j,\\
CF= F^\wedge(0),&\quad Dd= S^\wedge(0)d,\end{split}\end{equation} where $1\leq j\leq q,$ $F\in\mathbf{H}(S),$
$d\in\mathcal{C}_2.$ Then the following statements hold true:
\begin{enumerate}
\item The formulae \eqref{defabcd} define a bounded linear operator
\[V_j=\begin{pmatrix}A_j&B_j\\ C& D\end{pmatrix}\,:\,\begin{pmatrix}\mathbf{H}(S)\\\mathcal{C}_2\end{pmatrix}
\longrightarrow\begin{pmatrix}\mathbf{H}(S)\\\mathcal{C}_2\end{pmatrix},\] which satisfies
\begin{equation}\label{quaco}
V_jV_j^*=\begin{pmatrix}\hat{M}_{\alpha_j^*\alpha_j} &0\\ 0& I\end{pmatrix}.\end{equation} In particular, the
space $\mathbf{H}(S)$ is $A_j$-invariant for $j=1,\dots,q.$

\item The operators $A_j,B_j,C,D$ satisfy the relations
\begin{align}\label{belin}
A_\ell F=(A_jF)\alpha_j^*\alpha_\ell,&\quad B_\ell d=(B_jd)\alpha_j^*\alpha_\ell,\\
\label{melin}
A_j(Fc)=(A_jF)\alpha_j^*c\alpha_j,&\quad B_j(dc)=(B_jd)\alpha_j^*c\alpha_j,\\
C(Fc)= (CF)c,&\quad D(dc)= (Dd)c\label{endlin}\end{align} for every $F\in\mathbf{H}(S),$ $c\in\mathcal{C},$
$d\in\mathcal{C}_2,$ $1\leq j,\ell\leq q.$

\item
Let \begin{equation}\label{llg6}S=\sum_{w\in\mathcal{F}_q}w^*S_{[w]},\quad\text{where}\quad
S_{[w]}\in\mathcal{C}\end{equation} be the pointwise absolutely convergent expansion of $S$ as in Theorem
\ref{trivprop}. Then for every $d\in\mathcal{C}_2$ it holds that
\begin{equation}\label{coefrep}S_{[w]}{d}=\left\{\begin{aligned}
Dd,&\quad\text{if}\quad w=\oo,\\
 w(CA^v B_j d)v^*\alpha_j^*&\quad \forall w,v\in\mathcal{F}_q\,:\,|w|=|v|+1,\ \forall j\,:\,1\leq j\leq q,
 \end{aligned}\right.
\end{equation}
where
\begin{equation}\label{stnot8}
A^{\alpha_{i_1}\dots\alpha_{i_k}}\defi A_{i_k}\dots A_{i_1},\quad A^{\oo}\defi I.
\end{equation}
\end{enumerate}
\end{theorem}

\begin{proof}\ \paragraph{(I)}
Following  the idea  of \cite[Theorem 2.3]{adr1},   we define a linear operator \[\hat{V}_j:
\spa\left\{\begin{pmatrix}K_{S}^cd\\e\end{pmatrix}\,:\,c\in\mathbb{B}(\mathcal{T}),d,e\in\mathcal{C}_2\right\}
\longrightarrow\mathbf{H}(S)\oplus\mathcal{C}_2\] by
\begin{equation}\label{defhat}\hat{V}_j\begin{pmatrix}K_{S}^cd\\e\end{pmatrix}=\begin{pmatrix}\mathcal{B}_S\\CM_S^*\end{pmatrix}
(K_\wedge^cd\alpha_j^*+e).\end{equation}

 We claim that the operator $\hat{V}_j$ is well-defined and contractive; moreover,
 for every $
 F_1,F_2\in\spa\left\{\begin{pmatrix}K_{S}^cd\\e\end{pmatrix}\,:\,
 c\in\mathbb{B}(\mathcal{T}),d,e\in\mathcal{C}_2\right\}$
 \begin{equation}\label{relprop}
 \langle \hat{V}_jF_1,\hat{V}_jF_2\rangle_{\mathbf{H}(S)\oplus\mathcal{C}_2}=
 \left\langle\begin{pmatrix}\hat{M}_{\alpha_j^*\alpha_j} &0\\ 0&
 I\end{pmatrix}F_1,F_2\right\rangle_{\mathbf{H}(S)\oplus\mathcal{C}_2}.\end{equation}

Indeed, denote for the moment by  $C^*$ the adjoint of  $C$ in $\mathbf{H}_2(\mathcal{T})$ (that is,  the
injection of $\mathcal{C}_2$ into $\mathbf{H}_2(\mathcal{T})$).
 Then, in view of Definition \ref{debro} and the
statements (I) and (II) in Proposition \ref{glele}, we obtain
\begin{multline*}
\left\langle\hat{V}_j\begin{pmatrix}K_S^{c_1}d_1\\
e_1\end{pmatrix},\hat{V}_j\begin{pmatrix}K_S^{c_2}d_2\\
e_2\end{pmatrix}\right\rangle_{\mathcal{H}(S)\oplus\mathcal{C}_2}=\\
 =\left\langle\begin{pmatrix}\mathcal{B}_S\\CM_S^*\end{pmatrix}(K_\wedge^{c_1}d_1\alpha_j^*+e_1),
\begin{pmatrix}\mathcal{B}_S\\CM_S^*\end{pmatrix}(K_\wedge^{c_2}d_2\alpha_j^*+e_2)
\right\rangle_{\mathcal{H}(S)\oplus\mathcal{C}_2}\\= \langle(\mathcal{B}_S+M_S
C^*CM_S^*)(K_\wedge^{c_1}d_1\alpha_j^*+e_1), K_\wedge^{c_2}d_2\alpha_j^*+e_2\rangle_2\\= \langle(I-M_S
\hat{M}_{\alpha_j^*}A_jM_S^*)(K_\wedge^{c_1}d_1\alpha_j^*+e_1), K_\wedge^{c_2}d_2\alpha_j^*+e_2\rangle_2\\
=\langle(I-\hat{M}_{\alpha_j^*}M_S M_S^*A_j)(K_\wedge^{c_1}d_1\alpha_j^*+e_1),
K_\wedge^{c_2}d_2\alpha_j^*+e_2\rangle_2\\=\langle
A_j\hat{M}_{\alpha_j^*}K_\wedge^{c_1}d_1,K_\wedge^{c_2}d_2\rangle_2+\langle e_1,e_2\rangle_2- \langle A_j
\hat{M}_{\alpha_j^*}M_S M_S^*A_j\hat{M}_{\alpha_j^*}K_\wedge^{c_1}d_1,K_\wedge^{c_2}d_2\rangle_2\\=
\langle\hat{M}_{\alpha_j^*\alpha_j}\mathcal{B}_SK_\wedge^{c_1}d_1,K_\wedge^{c_2}d_2\rangle_2+\langle
e_1,e_2\rangle_2\\=\langle\hat{M}_{\alpha_j^*\alpha_j}K_S^{c_1}d_1,K_S^{c_2}d_2\rangle_{\mathbf{H}(S)}+\langle
e_1,e_2\rangle_2.
\end{multline*}

Thus \eqref{relprop} holds. Since, according  to  Proposition \ref{cla2}, the space $\mathbf{H}(S)$ is right
$\mathcal{C}$-invariant and $\|\alpha_j^*\alpha_j\|=1,$ we conclude that $\hat{V}_j$ is well-defined and
contractive. Since the span of $K_S^ck$ is dense in $\mathbf{H}(S)$ (see \eqref{dense2}), $\hat{V}_j$ can be
extended as a contraction on $\mathcal{H}(S)\oplus\mathcal{C}_2$, which satisfies
\[\hat{V}_j^*\hat{V}_j=\begin{pmatrix}\hat{M}_{\alpha_j^*\alpha_j} &0\\ 0&
 I\end{pmatrix}.\]

In order to complete the proof of statement (I), it suffices to observe that $V_j=\hat{V}_j^*.$

\paragraph{(II)}
The identities \eqref{belin} -- \eqref{endlin} follow immediately from \eqref{defabcd} and the Cuntz relations
\eqref{cuntz}.

\paragraph{(III)}
 The proof parallels the proof of the statement (III) in Proposition \ref{glele}. Let us fix a sequence of indices
$j=i_0,i_1,i_2,i_3,\dots,$ $1\leq i_n\leq q.$  Then for $d\in\mathcal{C}_2$ we have
\begin{multline*}Sd=Dd+(S-S_{[\oo]})d=Dd+(B_{j}d)\alpha_{j}^*=Dd+(CB_{j}d)\alpha_{j}^*+(A_{i_1}B_{j}d)
\alpha_{i_1}^*\alpha_{j}^*=\dots\\
=D{d}+\sum_{n=0}^m (CA^{\mu_n}B_j{d})\mu_n^*\alpha_j^*+(A^{\mu_{m+1}}B_j{d})\mu_{m+1}^*\alpha_j^*,\quad
 m=0,1,2,\dots\end{multline*}
 where  $\mu_0=\oo,$
$\mu_n=\alpha_{i_1}\dots\alpha_{i_n}.$ Now it follows from the identity \eqref{tayco} in Theorem \ref{trivprop},
applied to $Sd$, that
$$\sum_{|w|=n+1}w^*S_{[w]}d=(CA^{\mu_n}B_j{d})\mu_n^*\alpha_j^*\quad \forall n\geq 0$$
and hence
$$S_{[w]}d=w(CA^{\mu_{|w|-1}}B_j{d})\mu_{|w|}^*,\quad |w|\geq 1.$$
Since the sequence of indices $i_0,i_1,i_2,i_3,\dots$ was chosen arbitrarily, we obtain \eqref{coefrep}.
\end{proof}

A result which is converse to Theorem \ref{main1} can be formulated as follows:

\begin{theorem}\label{main2}
 Let $\mathbf{H}$ be a right $\mathcal{C}$-invariant Hilbert space
included in $\mathbf{H}_2(\mathcal{T})$. Assume that for some $j$, $1\leq j\leq q$, there exists a bounded
linear operator
\[V_j=\begin{pmatrix}A_j&B_j\\ C& D\end{pmatrix}\,:\,\begin{pmatrix}\mathbf{H}\\\mathcal{C}_2\end{pmatrix}
\longrightarrow\begin{pmatrix}\mathbf{H}\\\mathcal{C}_2\end{pmatrix},\] for which the relations \eqref{melin},
\eqref{endlin} and \eqref{quaco} hold true. Then: \begin{enumerate}
\item The series
\begin{equation}\label{finser}S=\sum_{w\in\mathcal{F}_q}w^*S_{[w]},\end{equation}
 where the coefficients $S_{[w]}\in\mathcal{C}$ are determined
by
\begin{equation}\label{fincoe}
\forall d\in\mathcal{C}_2\quad S_{[w]}d=\left\{\begin{aligned}Dd,&\quad\text{if}\quad w=\oo,\\
w(CA_j^{n-1}B_j d)\alpha_j^{n*},&\quad\text{if}\quad |w|=n\geq 1,\end{aligned}\right.\end{equation}
 defines (in the sense of pointwise absolute convergence) a Schur multiplier $S\in\mathbf{H}(\mathcal{T})$.

\item The series
\begin{equation}\label{finpoev}
E_cF=\sum_{n=0}^\infty(c\alpha)^n\cdot\left(CA_j^nF\right)\cdot \alpha_j^{n*},\quad
F\in\mathbf{H},c\in\mathbb{B}(\mathcal{T})\end{equation}    converges  absolutely in  $\mathcal{C}_2$ and
defines a bounded linear operator $E_c$ from $\mathbf{H}$ to $\mathcal{C}_2.$ This operator $E_c$ and the
 kernel $K_S^c$ defined in \eqref{repkers} satisfy
 \begin{equation}\label{repkerid}
 (K_S^d)^\wedge(c)k=E_cE_d^*k\quad\forall c,d\in\mathbb{B}(\mathcal{T}),k\in\mathcal{C}_2.\end{equation}

 \item If the operators $A_j$ and $C$ are as in \eqref{defabcd},  then $\mathbf{H}$ is the de Branges -- Rovnyak
 space associated with the Schur multiplier $S$:
 $$\mathbf{H}=\mathbf{H}(S).$$
 \end{enumerate}
\end{theorem}
\begin{proof}
Let  $c\in\mathbb{B}(\mathcal{T}).$ Then, as follows from \eqref{quaco} and \eqref{ball}, the series
\eqref{finpoev} is absolutely convergent in the $\mathcal{C}_2$-norm and defines a bounded linear operator $E_c$
from $\mathbf{H}$ to $\mathcal{C}_2$.\\

 Let us consider the linear operator
$s(c):\mathcal{C}_2\longrightarrow\mathcal{C}_2,$ defined by
\begin{equation}\label{fifi}
s(c)k\defi Dk+ \sum_{n=1}^\infty(c\alpha)^n\cdot(CA_j^{n-1}B_j k)\cdot\alpha_j^{n*},\quad
k\in\mathcal{C}_2.\end{equation}
 Here, in view of \eqref{quaco} and  \eqref{ball}, the series is absolutely
convergent in the $\mathcal{C}_2$-norm and the operator $s(c)$ is bounded. Moreover, as follows from
\eqref{melin} and \eqref{endlin}, the operator $s(c)$ commutes with $\hat{M}_d$ for every $d\in\mathcal{C}.$
Hence, according to the statement (II) of Lemma \ref{commut}, $s(c)\in\mathcal{C}.$ \\

Next we observe that
\[\begin{pmatrix} (M_{c\alpha}\hat{M}_{\alpha_j^*})E_c &
I\end{pmatrix}V_j=\begin{pmatrix}E_c&s(c)\end{pmatrix}.
\]
Note that, as follows from \eqref{melin} and \eqref{endlin}, the self-adjoint operator
$\hat{M}_{\alpha_j^*\alpha_j}$ commutes with $E_d$ for every $d\in\mathcal{C}.$ From \eqref{quaco} we obtain
\begin{multline*}
 \begin{pmatrix}E_c&s(c)\end{pmatrix}\begin{pmatrix}E_d&s(d)\end{pmatrix}^*=
\begin{pmatrix}(M_{c\alpha}\hat{M}_{\alpha_j^*})E_c & I\end{pmatrix}V_jV_j^*\begin{pmatrix}
(M_{d\alpha}\hat{M}_{\alpha_j^*})E_d &
I\end{pmatrix}^*\\=(M_{c\alpha}\hat{M}_{\alpha_j^*})E_c\hat{M}_{\alpha_j^*\alpha_j}E_d^*
(M_{\alpha^*d^*}\hat{M}_{\alpha_j})+I=(M_{c\alpha}\hat{M}_{\alpha_j^*})E_cE_d^*
(M_{\alpha^*d^*}\hat{M}_{\alpha_j})+I,\end{multline*} hence
\[I-s(c)s(d)^*=E_cE_d^*-
(M_{c\alpha}\hat{M}_{\alpha_j^*})E_cE_d^*(M_{\alpha^*d^*}\hat{M}_{\alpha_j}).\]

Therefore, for every $k\in\mathcal{C}_2$ the kernel $K_s(c,d),$ which appears in the equation \eqref{defks} of
Theorem \ref{Schthrm}, satisfies
\begin{multline*}K_s(c,d)k=\sum_{n=0}^\infty(c\alpha)^n(I-s(c)s(d)^*)(\alpha^*d^*)^nk\\=
\sum_{n=0}^\infty(c\alpha)^n\big((I-s(c)s(d)^*)((\alpha^*d^*)^nk\alpha_j^n)\big)\alpha_j^{n*}\\=
\sum_{n=0}^\infty(c\alpha)^n\big(E_cE_d^*(\alpha^*d^*)^nk\alpha_j^n\big)\alpha_j^{n*}-\\-
\sum_{n=0}^\infty(c\alpha)^{n+1}\big(E_cE_d^*(\alpha^*d^*)^{n+1}k\alpha_j^{n+1}\big)\alpha_j^{(n+1)*}=
 E_cE_d^*k.\end{multline*} In particular, the kernel $K_s(c,d)$ is positive. Now, according to Theorem \ref{Schthrm},
there exists a Schur multiplier $S\in\mathbf{H}(\mathcal{T})$ such that
\[s(c)=S^\wedge(c)\quad\forall c\in\mathbb{B}(\mathcal{T})\] and the identity \eqref{repkers}
 implies \eqref{repkerid}.\\

  Furthermore, in view of \eqref{fifi} and the
statement (II) of Lemma \ref{poevprop}, we obtain the formulae \eqref{fincoe} for the coefficients $S_{[w]}$ of
the pointwise absolutely convergent expansion
 \eqref{finser} of $S$. This completes the proof of the statements (I) and (II) of the Theorem.\\

In order to prove the statement (III), we note that if the operators $A_j$ and $C$ are as in \eqref{defabcd},
then, according to  the statement (III) of Proposition \ref{glele} and the definition \eqref{lepoev} of the
point evaluation,
\[E_cF=F^\wedge(c)\quad\forall F\in\mathbf{H}, c\in\mathbb{B}(\mathcal{T}.\]
Therefore, the identity \eqref{repkerid} implies \begin{align*}K_S^c k=E_c^*k\in\mathbf{H}&\quad \forall
c\in\mathbb{B}(\mathcal{T}),k\in\mathcal{C}_2;\\
\langle F, K_S^c k\rangle_{\mathbf{H}}=\langle F^\wedge(c),k\rangle_2,&\quad \forall F\in\mathbf{H},
c\in\mathbb{B}(\mathcal{T}),k\in\mathcal{C}_2.\end{align*} In view of Proposition \ref{cla2}  and the statement
(II) of Lemma \ref{poevprop},
$$\mathbf{H}=\overline{\spa}\{K_S^{c}{k}\,:\,{k}\in\mathcal{C}_2,c\in\mathbb{B}(\mathcal{T})\}
=\mathbf{H}(S)\quad\text{and}\quad \|\cdot\|_{\mathbf{H}}=\|\cdot\|_{\mathbf{H}(S)}.$$
\end{proof}

The formulae \eqref{coefrep} in Theorem \ref{main1} play the role of the backward shift realization
\eqref{real???} in the present setting: they allow to represent a given Schur multiplier $S$ as the transfer
operator of a multiscale input-state-output system, as described in the following theorem.

\begin{theorem}\label{isoform}
Let \(S\in\mathbf{H}(\mathcal{T})\) be a Schur multiplier and let $\mathbf{H}(S)$ be the associated de Branges
-- Rovnyak space. Let  the operators \[V_j=\begin{pmatrix}A_j&B_j\\
C& D\end{pmatrix}\,:\,\begin{pmatrix}\mathbf{H}(S)\\\mathcal{C}_2\end{pmatrix}
\longrightarrow\begin{pmatrix}\mathbf{H}(S)\\\mathcal{C}_2\end{pmatrix},\quad 1\leq j\leq q,\] be defined by
\eqref{defabcd} as in Theorem \ref{main1}. Let $U\in\mathbf{H}_2(\mathcal{T}),$
\[U=\sum_{w\in\mathcal{F}_q} w^* U_{[w]},\quad\text{where}\quad U_{[w]}\in\mathcal{C}_2,\]
and let $Y=SU\in\mathbf{H}_2(\mathcal{T})$. Then the coefficients $Y_{[w]}\in\mathcal{C}_2$
 of the expansion
\[ Y=\sum_{w\in\mathcal{F}_q} w^*Y_{[w]}\]
satisfy the recurrent relations
\begin{equation}\label{isoeq}\left\{\begin{aligned}
X_\oo&=0,\\
X_{w\alpha_j}&=A_jX_w+B_jU_w,\\
Y_{[w]}&=w(CX_v+DU_v)v^*\quad \forall v\in\mathcal{F}_q\,:\,|v|=|w|,
\end{aligned}\right.\end{equation}
where
\[U_w\defi\sum_{v:\,|v|=|w|}v^*U_{[v]}w.\]
\end{theorem}
\begin{proof}
The case $w=\oo$ is trivial, so let us assume $n=|w|\geq 1.$
 Let $v=\alpha_{i_1}\dots\alpha_{i_n}$ and
denote
$$ v_{k}=\alpha_{i_1}\dots\alpha_{i_k},\quad v_{\overline{k}}=\alpha_{i_{k+1}}\dots\alpha_{i_n},\quad
v_{0}=v_{\overline{n}}=\oo.$$

As follows from \eqref{coefrep},
\begin{multline*}
\sum_{\substack{u\in\mathcal{F}_q\\|u|=n}}u^*Y_{[u]}=\sum_{\substack{\mu,\nu\in\mathcal{F}_q\\|\mu|+|\nu|=n}}\mu^*
S_{[\mu]}\nu^*U_{[\nu]}\\=\sum_{k=0}^{n-1}\sum_{\substack{u\in\mathcal{F}_q\\|u|=k}}CA^{v_{\overline{k+1}}}
B_{i_{k+1}}(u^*U_{[u]} v_{k})v^*+\sum_{\substack{u\in\mathcal{F}_q\\|u|=n}} D(u^*U_{[u]}
v)v^*\\=\left(C\sum_{k=0}^{n-1} A^{v_{\overline{k+1}}}B_{i_{k+1}}U_{v_{k}}+DU_v\right)v^*,
\end{multline*}
hence, according to \eqref{cuntz},
\[Y_{[w]}=w\left(C\sum_{0\leq k\leq n-1}
A^{v_{\overline{k+1}}}B_{i_{k+1}}U_{v_{k}}+DU_v\right)v^*.\] \ \\

 Denote \[X_v=\sum_{0\leq k\leq
|v|-1}A^{v_{\overline{k+1}}}
B_{i_{k+1}}U_{v_{k}},\quad v\in\mathcal{F}_q,|v|\geq 1.\] Then
\[
X_{v\alpha_j}=\sum_{0\leq k\leq
|v|-1}A_jA^{v_{\overline{k+1}}}B_{i_{k+1}}U_{v_{k}}+B_jU_{v}=A_jX_v+B_jU_v.\]
In the case $v=\alpha_j$ we
have \[X_{\alpha_j}=B_jU_{\oo}=A_jX_{\oo}+B_jU_{\oo}.\] Thus the relations \eqref{isoeq} hold.
\end{proof}

\section{The Blaschke factors}\label{s7}
In this section we present an important example of Schur multiplier, which plays a role in interpolation. We
follow the ideas of \cite[p. 86--90]{MR93b:47027}. Let \({c}\in {\mathbb B}({\mathcal T})\) and consider the
operator
\[R_c\defi\sum_{n=0}^\infty (c\alpha )^n(\alpha^*
c^*)^n=(K_\wedge^c)^\wedge(c)\in\mathcal{C}.\] Then $R_c>0$ and,
moreover,\begin{equation}\label{imin}R_c=I+c(\alpha R_c\alpha^*)c^*.\end{equation} Hence $R_c>I$. Next we define
\[L_c\defi\alpha(R_c-R_c\alpha^* c^*R_c^{-1}c\alpha
R_c)\alpha^*\in\mathcal{C}^{q\times q}.\]
\begin{proposition}
The following holds:
\begin{align}
L_c&>0,\\
L_c^{-1}&=c^*c+\alpha R_c^{-1}\alpha^*=I+c^*c-\alpha c L_c c^*\alpha^*,\\
cL_c&=R_c^{-1}c\alpha R_c\alpha^*.
\end{align}
\end{proposition}
\begin{proof}
Taking into account \eqref{imin}, \begin{multline*}L_c(c^*c+\alpha R_c^{-1}\alpha^*)=\alpha R_c\alpha^*
c^*c+I-\alpha R_c\alpha^* c^*R_c^{-1}c\alpha R_c\alpha^* c^*c-\alpha R_c\alpha^* c^*R_c^{-1}c\\=\alpha
R_c\alpha^* c^*c+I-\alpha R_c\alpha^* c^*R_c^{-1}(R_c-I)c-\alpha R_c\alpha^* c^*R_c^{-1}c=I.\end{multline*}
Since $c^*c+\alpha R_c^{-1}\alpha^*$ is positive definite, this means
$$L_c=(c^*c+\alpha R_c^{-1}\alpha^*)^{-1}>0.$$
The rest of the identities follow from \eqref{imin} analogously.
\end{proof}
Note that
\[(L_cc^*\alpha^*)^n=(\alpha R_c\alpha^* c^*
R_c^{-1}\alpha^*)^n=\alpha R_c\alpha^* c^* (\alpha^* c^*)^{n-1}R_c^{-1}\alpha^*,\] hence \(\alpha L_c
c^*\alpha^*\in\mathbb{B}(\mathcal{T})\).

\begin{definition}
Let \({c}\in {\mathbb B}({\mathcal T})\).
 The operator
\begin{equation*}
B_{c}= (\alpha^*-c)(1-L_{c}c^*\alpha^*)^{-1}\sqrt{L_c}\in{\bf H}(\mathcal T)
\end{equation*}
is called the  Blaschke factor, corresponding to the constant $c$.
\end{definition}

\begin{proposition}\label{blaprop}
The operator $B_c$ is unitary. In particular, the multiplication operator ${\mathcal
M}_{B_{c}}:\mathbf{H}_2(\mathcal{T})^q\longrightarrow\mathbf{H}_2(\mathcal{T})$ is an isometry.
\end{proposition}
\begin{proof}
We have
\[(I-\alpha c L_c)L_c^{-1}(I-L_c
c^*\alpha^*)=(\alpha-c^*)(\alpha^*-c),\] hence $B_c^*B_c=I$ and $M_{B_c}$ is an isometry. Furthermore,
\[(\alpha^*-c)=\alpha^*(I-\alpha c)\] has a bounded inverse (not
causal), hence $B_c$ is also invertible and unitary.
\end{proof}
\begin{theorem}
Let $F\in\mathbf{H}_2(\mathcal{T})$, $c\in\mathbb{B}(\mathcal{T}).$ Then $F^\wedge(c)=0$ if and only if $F$ is
of the form
$$F=B_c\cdot G,$$
where $G\in\mathbf{H}_2(\mathcal{T})^q,$ $\|G\|_2=\|F\|_2$.
\end{theorem}
\begin{proof}
First, let us assume that $F=B_c\cdot G,$ where $G\in\mathbf{H}_2(\mathcal{T})^q.$ Then $\|F\|_2=\|G\|_2$ by
Proposition \ref{blaprop}. Furthermore, as follows from Lemma
\ref{poevprop}, $B_c^\wedge(c)=0$ and $F^\wedge(c)=(B_cG)^\wedge(c)=0$.\\

Conversely, assume that $F$ has the expansion $F=\sum_{w\in\mathcal{F}_q}w^*F_{[w]}$ and that $F^\wedge(c)=0$.
Then $F$ is represented by the \ series
$$F=F-F^\wedge(c)=\sum_{w\in\mathcal{F}_q}(w^*-(w^*)^\wedge(c))F_{[w]}.$$
Denote $G^\prime\defi(I-\alpha c)^{-1}\alpha F.$ Then $F=(\alpha^*-c)G^\prime,$
$G\in\mathbf{X}_2(\mathcal{T})^q$ and, moreover,
\begin{multline*}G^\prime=\sum_{w\in\mathcal{F}_q}(I-\alpha
c)^{-1}\alpha(w^*-(w^*)^\wedge(c))F_{[w]}\\= \sum_{\substack{w\in\mathcal{F}_q\\|w|\geq
1}}\alpha\left(I+c\alpha+(c\alpha)^2+\dots+(c\alpha)^{|w|-1}\right)w^*F_{[w]}.\end{multline*} Since each term of
the last series is causal, $G^\prime\in\mathbf{H}_2(\mathcal{T})^q.$ It remains to define
$$G\defi\sqrt{L_c^{-1}}(I-L_c c^*\alpha^*)G^\prime$$
to complete the proof.
\end{proof}
\bibliographystyle{plain}


\end{document}